\documentclass[11pt]{article}
\usepackage{graphicx}
\usepackage{amsmath}
\usepackage{amssymb}
\usepackage{theorem}
\usepackage{euscript}
\usepackage{epic,eepic}
\usepackage{pstricks}
\title{\sffamily A PARALLEL SPLITTING METHOD FOR\\
WEAKLY COUPLED MONOTONE INCLUSIONS\footnote{Contact author: 
P. L. Combettes, {\ttfamily plc@math.jussieu.fr},
phone: +33 1 4427 6319, fax: +33 1 4427 7200.
This work was supported by the Agence Nationale de la Recherche under 
grant ANR-08-BLAN-0294-02.}}
\author{H\'edy Attouch,$^1$ Luis M. Brice\~{n}o-Arias,$^2$ and 
Patrick L. Combettes$^3$
\\[5mm]
\small
$\!^1$Universit\'e Montpellier II\\
\small Institut de Math\'ematiques et de Mod\'elisation de 
Montpellier -- UMR 5149\\
\small 34095 Montpellier Cedex 5, France 
(attouch@math.univ-montp2.fr)\\[4mm]
\small $\!^2$UPMC Universit\'e Paris 06\\
\small \'Equipe Combinatoire et Optimisation -- UMR 7090\\
\small 75005 Paris, France (lbriceno@math.jussieu.fr)\\[4mm]
\small $\!^3$UPMC Universit\'e Paris 06\\
\small Laboratoire Jacques-Louis Lions -- UMR 7598\\
\small 75005 Paris, France (plc@math.jussieu.fr)
}

\date{\ttfamily ~}
\newcommand{\Frac}[2]{\displaystyle{\frac{#1}{#2}}} 
\newcommand{\Scal}[2]{{\bigg\langle{{#1}\:\bigg |~{#2}}\bigg\rangle}}
\newcommand{\scal}[2]{{\left\langle{{#1}\mid{#2}}\right\rangle}}
\newcommand{\pscal}[2]{\langle\langle{#1}\mid{#2}\rangle\rangle} 
\newcommand{\Menge}[2]{\Big\{{#1}~\big |~{#2}\Big\}} 
\newcommand{\menge}[2]{\big\{{#1}~\big |~{#2}\big\}} 
\newcommand{\moyo}[2]{\ensuremath{\sideset{^{#2}}{}{\operatorname{}}\!\!#1}}

\newcommand{\HH}{\ensuremath{{\mathcal H}}}
\newcommand{\HHH}{\ensuremath{\boldsymbol{\mathcal H}}}

\newcommand{\GG}{\ensuremath{{\mathcal G}}}

\newcommand{\BL}{\ensuremath{\EuScript B}\,}

\newcommand{\emp}{\ensuremath{{\varnothing}}}

\newcommand{\Id}{\ensuremath{\operatorname{Id}}\,}
\newcommand{\RR}{\ensuremath{\mathbb{R}}}

\newcommand{\haus}{\ensuremath{\operatorname{haus}}}
\newcommand{\bdry}{\ensuremath{\operatorname{bdry}}}
\newcommand{\RP}{\ensuremath{\left[0,+\infty\right[}}

\newcommand{\RPP}{\ensuremath{\left]0,+\infty\right[}}

\newcommand{\RXX}{\ensuremath{\left[-\infty,+\infty\right]}}
\newcommand{\RX}{\ensuremath{\left]-\infty,+\infty\right]}}
\newcommand{\NN}{\ensuremath{\mathbb N}}

\newcommand{\KK}{\ensuremath{\mathbb K}}

\newcommand{\cart}{\ensuremath{\raisebox{-0.5mm}{\mbox{\LARGE{$\times$}}}}\!}

\newcommand{\gr}{\ensuremath{\operatorname{gra}}}
\newcommand{\exi}{\ensuremath{\exists\,}}
\newcommand{\pinf}{\ensuremath{{+\infty}}}
\newcommand{\minf}{\ensuremath{{-\infty}}}
\newcommand{\dom}{\ensuremath{\operatorname{dom}}}
\newcommand{\prox}{\ensuremath{\operatorname{prox}}}

\newcommand{\inte}{\ensuremath{\operatorname{int}}}

\newtheorem{theorem}{Theorem}[section]
\newtheorem{lemma}[theorem]{Lemma}
\newtheorem{corollary}[theorem]{Corollary}
\newtheorem{proposition}[theorem]{Proposition}
\theoremstyle{plain}{\theorembodyfont{\rmfamily}%
}
\theoremstyle{plain}{\theorembodyfont{\rmfamily}%
\newtheorem{algorithm}[theorem]{Algorithm}}
\theoremstyle{plain}{\theorembodyfont{\rmfamily}%
\newtheorem{example}[theorem]{Example}}
\theoremstyle{plain}{\theorembodyfont{\rmfamily}%
\newtheorem{remark}[theorem]{Remark}}
\theoremstyle{plain}{\theorembodyfont{\rmfamily}%
\newtheorem{definition}[theorem]{Definition}}
\theoremstyle{plain}{\theorembodyfont{\rmfamily}%
\newtheorem{problem}[theorem]{Problem}}

\numberwithin{equation}{section}
\begin{document}
\maketitle

\begin{abstract}
A parallel splitting method is proposed for solving systems of coupled 
monotone inclusions in Hilbert spaces. Convergence is established for 
a wide class of coupling schemes. Unlike classical alternating 
algorithms, which are limited to two variables and linear coupling, our 
parallel method can handle an arbitrary number of variables as well as 
nonlinear coupling schemes. The breadth and flexibility of the proposed 
framework is illustrated through applications in the areas of evolution 
inclusions, dynamical games, signal recovery, image decomposition, best
approximation, network flows, and variational problems in Sobolev 
spaces.
\end{abstract}


\newpage
\section{Problem statement}
\label{sec:1}

This paper is concerned with the numerical solution of systems of
coupled monotone inclusions in Hilbert spaces. A simple instance 
of this problem is to
\begin{equation}
\label{e:prob63}
\text{find}\;x_1\in\HH,\, x_2\in\HH\quad\text{such that}\quad
\begin{cases}
0\in A_1x_1+x_1-x_2\\
0\in A_2x_2+x_2-x_1,
\end{cases}
\end{equation}
where $(\HH,\|\cdot\|)$ is a real Hilbert space, and where $A_1$ and 
$A_2$ are maximal monotone operators acting on $\HH$. This 
formulation arises in various areas of nonlinear analysis 
\cite{Nona05}. For example, if $A_1=\partial f_1$ and 
$A_2=\partial f_2$ are the subdifferentials of proper lower 
semicontinuous convex functions $f_1$ and $f_2$ from $\HH$ to 
$\RX$, \eqref{e:prob63} is equivalent to 
\begin{equation}
\label{e:prob62}
\underset{x_1\in\HH,\,x_2\in\HH}{\mathrm{minimize}}\;\;
f_1(x_1)+f_2(x_2)+\frac{1}{2}\|x_1-x_2\|^2.
\end{equation}
This joint minimization problem, which was first investigated in 
\cite{Acke80}, models problems in disciplines such as the 
cognitive sciences \cite{Atto07}, 
image processing \cite{Smms05}, and signal processing \cite{Gold85} (see 
also the references therein for further applications in mechanics, 
filter design, and dynamical games). In particular, if $f_1$ and $f_2$ are
the indicator functions of closed convex subsets $C_1$ and $C_2$ of 
$\HH$, \eqref{e:prob62} reverts to the classical best 
approximation pair problem \cite{Baus04,Che59a,Gubi67}
\begin{equation}
\label{e:prob61}
\underset{x_1\in C_1,\,x_2\in C_2}{\mathrm{minimize}}\;\;
\|x_1-x_2\|.
\end{equation}
On the numerical side, a simple 
algorithm is available to solve \eqref{e:prob63}, namely,
\begin{equation}
\label{e:nona05}
x_{1,0}\in\HH
\quad\text{and}\quad 
(\forall n\in\NN)\quad 
\begin{cases}
x_{2,n}&\!\!\!\!=(\Id+{A_2})^{-1}x_{1,n}\\
x_{1,n+1}&\!\!\!\!=(\Id+{A_1})^{-1}x_{2,n}.
\end{cases}
\end{equation}
This alternating resolvent method produces sequences 
$(x_{1,n})_{n\in\NN}$ and $(x_{2,n})_{n\in\NN}$ that converge weakly 
to points $x_1$ and $x_2$, respectively, such that $(x_1,x_2)$ solves 
\eqref{e:prob63} if solutions exist \cite[Theorem~3.3]{Nona05}. In 
\cite{Atto08}, the variational formulation \eqref{e:prob62} was 
extended to
\begin{equation}
\label{e:prob5-}
\underset{x_1\in\HH_1,\,x_2\in\HH_2}{\mathrm{minimize}}\;\;
f_1(x_1)+f_2(x_2)+\frac12\|L_{1}x_1-L_{2}x_2\|^2_{\GG},
\end{equation}
where $\HH_1$, $\HH_2$, and $\GG$ are Hilbert spaces, 
$f_1\colon\HH_1\to\RX$ and $f_2\colon\HH_2\to\RX$ are proper 
lower semicontinuous convex functions, and $L_1\colon\HH_1\to\GG$ 
and $L_2\colon\HH_2\to\GG$ are linear and bounded. This problem
was solved in
\cite{Atto08} via an inertial alternating minimization procedure first
proposed in \cite{Atto07} for the case of the strongly coupled problem
\eqref{e:prob62}.

The above problems and their solution algorithms are limited to two 
variables which, in addition, must be linearly coupled. These are
serious restrictions since models featuring more than two variables 
and/or nonlinear coupling schemes arise naturally in applications.
The purpose of this paper is to address simultaneously these
restrictions by proposing a parallel algorithm for solving systems of
monotone inclusions involving an arbitrary number of variables and 
nonlinear coupling. The 
breadth and flexibility of this framework will be illustrated through 
applications in the areas of evolution inclusions, dynamical games, 
signal recovery, image decomposition, best approximation, 
network flows, and decomposition methods in Sobolev spaces.

We now state our problem formulation and our standing assumptions.

\begin{problem}
\label{prob:64}
Let $(\HH_i)_{1\leq i\leq m}$ be real Hilbert spaces, where $m\geq 2$. 
For every $i\in\{1,\ldots,m\}$, let 
$A_i\colon\HH_i\to 2^{\HH_i}$ be maximal monotone and let 
$B_{i}\colon\HH_1\times\cdots\times\HH_m\to\HH_i$. It is assumed that 
there exists $\beta\in\RPP$ such that
\begin{multline}
\label{e:genna07-21}
(\forall (x_1,\ldots,x_m)\in\HH_1\times\cdots\times\HH_m)
(\forall (y_1,\ldots,y_m)\in\HH_1\times\cdots\times\HH_m)\\
\sum_{i=1}^m\scal{B_{i}(x_1,\ldots,x_m)-B_{i}(y_1,\ldots,y_m)}
{x_i-y_i}\geq\beta\sum_{i=1}^m\big\|B_{i}(x_1,\ldots,x_m)-
B_{i}(y_1,\ldots,y_m)\big\|^2.
\end{multline}
The problem is to 
\begin{equation}
\label{e:genna07-1}
\text{find}\;x_1\in\HH_1,\ldots,x_m\in\HH_m\quad\text{such that}\quad
\begin{cases}
0\in A_1x_1+B_1(x_1,\ldots,x_m)\\
\quad\vdots\\
0\in A_mx_m+B_m(x_1,\ldots,x_m),
\end{cases}
\end{equation}
under the assumption that such points exist.
\end{problem}

In abstract terms, the system of inclusions in \eqref{e:genna07-1} 
models an equilibrium involving $m$ variables in different Hilbert 
spaces. The $i$th inclusion in this system is a perturbation of 
the basic inclusion $0\in A_ix_i$ by addition of the coupling term 
$B_{i}(x_1,\ldots,x_m)$. This type of coupling will be referred to 
as \emph{weak} in that it is not restricted to a simple linear 
combination of the variables as in \eqref{e:prob63}. 
As will be seen in Section~\ref{sec:3},
our analysis captures various linear and nonlinear coupling schemes. 
For example, if
\begin{equation}
(\forall i\in\{1,\ldots,m\})\quad\HH_i=\HH\quad\text{and}\quad
(\forall x\in\HH)\quad B_{i}(x,\ldots,x)=0,
\end{equation}
then Problem~\ref{prob:64} is a relaxation of the standard problem 
\cite{Cras95,Lema99} of finding a common zero of the operators 
$(A_i)_{1\leq i\leq m}$, i.e., of solving the inclusion 
$0\in\bigcap_{i=1}^mA_ix$. In particular, if $m=2$, $\HH_1=\HH_2=\HH$, 
$B_{1}=-B_{2}\colon(x_1,x_2)\mapsto x_1-x_2$, and $\beta=1/2$, then
Problem~\ref{prob:64} reverts to \eqref{e:prob63}. On the other hand, 
if
$m=2$, $A_1=\partial f_1$, $A_2=\partial f_2$, 
$B_1\colon(x_1,x_2)\mapsto L_1^*(L_1x_1-L_2x_2)$,  
$B_2\colon(x_1,x_2)\mapsto -L_2^*(L_1x_1-L_2x_2)$, and 
$\beta=(\|L_1\|^2+\|L_2\|^2)^{-1}$, then Problem~\ref{prob:64} 
reverts to \eqref{e:prob5-}.

The paper is organized as follows. In Section~\ref{sec:2}, we present
our algorithm for solving Problem~\ref{prob:64} and prove its 
convergence to solutions to Problem~\ref{prob:64}. 
In Section~\ref{sec:2+}, we describe various instances of 
\eqref{e:genna07-1} resulting from specific choices for the
operators $(A_i)_{1\leq i\leq m}$, e.g., minimization problems, 
variational inequalities, saddle-point problems, and evolution 
inclusions. 
In Section~\ref{sec:3}, we discuss examples of linear and nonlinear 
coupling schemes that can be obtained through specific choices of 
the operators $(B_i)_{1\leq i\leq m}$
in Problem~\ref{prob:64}. Applications to systems of evolution inclusions
are treated in Section~\ref{sec:6}. Section~\ref{sec:5} is devoted to 
variational formulations deriving from Problem~\ref{prob:64} and 
features various special cases. The applications treated in that section
include dynamical games, signal recovery, image decomposition, best
approximation, and network flows. Finally, Section~\ref{sec:7} 
describes an application to decomposition methods in Sobolev spaces. 

{\bfseries Notation.}
Throughout, $\HH$ and $(\HH_i)_{1\leq i\leq m}$ are real Hilbert spaces. 
Their scalar products are denoted by $\scal{\cdot}{\cdot}$ and the 
associated norms by $\|\cdot\|$. Moreover, $\Id$ denotes the identity 
operator on these spaces. The indicator function of a subset $C$ of 
$\HH$ is
\begin{equation}
\label{e:iota}
\iota_C\colon x\mapsto
\begin{cases}
0,&\text{if}\;\;x\in C;\\
\pinf,&\text{if}\;\;x\notin C
\end{cases}
\end{equation}
and the distance from $x\in\HH$ to $C$ is $d_C(x)=\inf_{y\in C}\|x-y\|$; 
if $C$ is nonempty closed and convex, the projection of $x$ onto $C$ is 
the unique point $P_Cx$ in $C$ such that $\|x-P_Cx\|=d_C(x)$. We denote 
by $\Gamma_0(\HH)$ the class of lower semicontinuous convex functions 
$f\colon\HH\to\RX$ which are proper in the sense that
$\dom f=\menge{x\in\HH}{f(x)<\pinf}\neq\emp$.
The subdifferential of $f\in\Gamma_0(\HH)$ is the maximal monotone 
operator
\begin{equation}
\label{e:subdiff}
\partial f\colon\HH\to 2^{\HH}\colon x\mapsto
\menge{u\in\HH}{(\forall y\in\HH)\;\;\scal{y-x}{u}+f(x)\leq f(y)}.
\end{equation}
If $\GG$ is a real Hilbert space, $\BL(\HH,\GG)$ is the space of 
bounded linear operators from $\HH$ to $\GG$ and $\BL(\HH)=\BL(\HH,\HH)$.
We denote by $\gr A=\menge{(x,u)\in\HH\times\HH}{u\in Ax}$ the graph
of a set-valued operator $A\colon\HH\to 2^{\HH}$, by
$\dom A=\menge{x\in\HH}{Ax\neq\emp}$ its domain, and by
$J_A=(\Id+A)^{-1}$ its resolvent. If $A$ is monotone, then $J_A$ is 
single-valued and nonexpansive and, furthermore, if $A$ is maximal
monotone, then $\dom J_A=\HH$. For complements and 
further background on convex analysis and monotone operator theory, 
see \cite{Aubi90,Brez73,Simo08,Zali02,Zeid90}.

\section{Algorithm}
\label{sec:2}
Let us start with a characterization of the solutions to
Problem~\ref{prob:64}.

\begin{proposition}
\label{p:10octobre2007}
Let $(x_i)_{1\leq i\leq m}\in\HH_1\times\cdots\times\HH_m$,
let $(\lambda_i)_{1\leq i\leq m}\in\left[0,1\right[^m$, and
let $\gamma\in\RPP$. 
Then $(x_i)_{1\leq i\leq m}$ solves
Problem~\ref{prob:64} if and only if
\begin{equation}
\label{e:10octobre2007}
(\forall i\in\{1,\ldots,m\})\quad
x_i=\lambda_ix_i+(1-\lambda_i)J_{\gamma A_i}\big(x_i-
\gamma B_{i}(x_1,\ldots,x_m)\big).
\end{equation}
\end{proposition}
\begin{proof}
Let $i\in\{1,\ldots,m\}$. Then, since $B_i$ is single-valued,
\begin{eqnarray}
0\in A_ix_i+B_{i}(x_1,\ldots,x_m)
&\Leftrightarrow&
x_i-\gamma B_{i}(x_1,\ldots,x_m)\in x_i+\gamma A_ix_i\nonumber\\
&\Leftrightarrow&
x_i=J_{\gamma A_i}\big(x_i-\gamma B_{i}(x_1,\ldots,x_m)\big)\nonumber\\
&\Leftrightarrow&
x_i=x_i+(1-\lambda_i)\big(J_{\gamma A_i}\big(x_i-\gamma B_{i}
(x_1,\ldots,x_m)\big)-x_i\big),
\end{eqnarray}
and we obtain \eqref{e:10octobre2007}.
\end{proof}

The above characterization suggests the following algorithm, which
constructs $m$ sequences $((x_{i,n})_{n\in\NN})_{1\leq i\leq m}$.
Recall that $\beta$ is the constant appearing in \eqref{e:genna07-21}.
\begin{algorithm}
\label{algo:1}
Fix $\varepsilon\in\left]0,\min\{1,\beta\}\right[$,
$(\gamma_n)_{n\in\NN}$ in $\left[\varepsilon,2\beta-\varepsilon\right]$,
$(\lambda_n)_{n\in\NN}$ in $\left[0,1-\varepsilon\right]$, and
$(x_{i,0})_{1\leq i\leq m}$ $\in\HH_1\times\cdots\times\HH_m$.
Set, for every $n\in\NN$,
\begin{equation}
\label{e:genna07-2}
\begin{cases}
x_{1,n+1}\!\!\!\!\!&=\lambda_{1,n}x_{1,n}
+(1-\lambda_{1,n})\Big(J_{\gamma_{n} A_{1,n}}\big(
x_{1,n}-\gamma_{n}(B_{1,n}(x_{1,n},\ldots,x_{m,n})+b_{1,n})\big)
+a_{1,n}\Big)\\
&~\vdots\\
x_{m,n+1}\!\!\!\!\!&=\lambda_{m,n}x_{m,n}
+(1-\lambda_{m,n})\Big(J_{\gamma_{n} A_{m,n}}\big(
x_{m,n}-\gamma_{n}(B_{m,n}(x_{1,n},\ldots,x_{m,n})+b_{m,n})\big)
+a_{m,n}\Big),
\end{cases}
\end{equation}
where, for every $i\in\{1,\ldots,m\}$, the following hold.
\begin{enumerate}
\item
\label{c:1}
$(A_{i,n})_{n\in\NN}$ are maximal monotone operators from 
$\HH_i$ to $2^{\HH_i}$ such that
\begin{equation}
\label{e:lc1}
(\forall\rho\in\RPP)\;\;\sum_{n\in\NN}\:\sup_{\|y\|\leq\rho}
\|J_{\gamma_n A_{i,n}}y-J_{\gamma_nA_i}y\|<\pinf.
\end{equation}
\item
\label{c:2}
$(B_{i,n})_{n\in\NN}$ are operators from 
$\HH_1\times\cdots\times\HH_m$ to $\HH_i$ such that
\begin{enumerate}
\item
\label{c:21}
the operators $(B_{i,n}-B_i)_{n\in\NN}$ are Lipschitz-continuous with 
respective constants $(\kappa_{i,n})_{n\in\NN}$ in $\RPP$ satisfying
$\sum_{n\in\NN}\kappa_{i,n}<\pinf$; and 
\item
\label{c:22}
there exists
${\boldsymbol z}\in\HH_1\times\cdots\times\HH_m$, independent of $i$,
such that $(\forall n\in\NN)$ 
$B_{i,n}{\boldsymbol z}=B_i{\boldsymbol z}$.
\end{enumerate}
\item
\label{c:3}
$(a_{i,n})_{n\in\NN}$ and $(b_{i,n})_{n\in\NN}$ are sequences in 
$\HH_i$ such that $\sum_{n\in\NN}\|a_{i,n}\|<\pinf$ and 
$\sum_{n\in\NN}\|b_{i,n}\|<\pinf$.
\item
\label{c:4}
$(\lambda_{i,n})_{n\in\NN}$ is a sequence in $\left[0,1\right[$
such that $\sum_{n\in\NN}|\lambda_{i,n}-\lambda_n|<\pinf$.
\end{enumerate}
\end{algorithm}

Conditions~\ref{c:1} and \ref{c:2} describe the types of approximations
to the original operators $(A_i)_{1\leq i\leq m}$ and 
$(B_i)_{1\leq i\leq m}$ which can be utilized.
Examples of approximations will be provided in 
Proposition~\ref{p:2008-07-09} and Remark~\ref{r:82}, respectively.
Condition~\ref{c:3} quantifies the tolerance which is
allowed in the implementation of these approximations 
(see \cite{Sico04,Hehe03,Kapl00} for specific 
examples), while \ref{c:4} quantifies that allowed in the
agent-dependent departure from the global relaxation scheme.
The parallel nature of Algorithm~\ref{algo:1} stems from the fact that
the $m$ evaluations of the resolvent operators in \eqref{e:genna07-2} 
can be performed independently and, therefore, simultaneously on
concurrent processors.

Our asymptotic analysis of Algorithm~\ref{algo:1} requires the 
following result on the convergence of the forward-backward algorithm.
This algorithm finds its roots in the projected gradient method 
\cite{Levi66} and certain methods for solving variational inequalities
\cite{Ausl72,Brez67,Merc79,Sibo70} (see also the bibliography of 
\cite{Opti04} for more recent developments). First, we need to define 
the notion of cocoercivity.

\begin{definition}
\label{d:coco}
Let $\chi\in\RPP$. An operator $B\colon\HH\to\HH$ is $\chi$-cocoercive if
\begin{equation}
\label{e:cocoercive}
(\forall x\in\HH)(\forall y\in\HH)\quad
\scal{x-y}{Bx-By}\geq\chi\|Bx-By\|^2.
\end{equation}
If $\chi=1$ in \eqref{e:cocoercive}, then $B$ is firmly nonexpansive. 
\end{definition}

\begin{lemma}{\rm\cite[Corollary~6.5]{Opti04}}
\label{l:10}
Let $(\HHH,|||\cdot|||)$ be a real Hilbert space, let $\chi\in\RPP$, let
${\boldsymbol A}\colon\HHH\to 2^{\HHH}$ be a maximal monotone operator, 
and let ${\boldsymbol B}\colon\HHH\to\HHH$ be a $\chi$-cocoercive 
operator such that 
$({\boldsymbol A}+{\boldsymbol B})^{-1}({\boldsymbol 0})\neq\emp$.
Fix $\varepsilon\in\left]0,\min\{1,\chi\}\right[$\,,
let $(\gamma_n)_{n\in\NN}$ be a sequence in 
$[\varepsilon,2\chi-\varepsilon]$, 
let $(\lambda_n)_{n\in\NN}$ be a sequence in 
$\left[0,1-\varepsilon\right]$,
and let $({\boldsymbol a}_n)_{n\in\NN}$ and 
$({\boldsymbol b}_n)_{n\in\NN}$ 
be sequences in $\HHH$ such that
$\sum_{n\in\NN}|||{\boldsymbol a}_n|||<\pinf$ and
$\sum_{n\in\NN}|||{\boldsymbol b}_n|||<\pinf$. Fix 
${\boldsymbol x}_0\in\HHH$ and, for every $n\in\NN$, set
\begin{equation}
\label{e:main2}
{\boldsymbol x}_{n+1}=\lambda_n{\boldsymbol x}_n+(1-\lambda_n)
\big(J_{\gamma_n {\boldsymbol A}}({\boldsymbol x}_n
-\gamma_n({\boldsymbol B}x_n+{\boldsymbol b}_n))+
{\boldsymbol a}_n\big).
\end{equation}
Then $({\boldsymbol x}_n)_{n\in\NN}$ converges weakly to a point in
$({\boldsymbol A}+{\boldsymbol B})^{-1}({\boldsymbol 0})$.
\end{lemma}

We shall also use the following fact.
\begin{lemma}{\rm\cite[Lemma~2.3]{Opti04}}
\label{l:9}
Let $(\HHH,|||\cdot|||)$ be a real Hilbert space, let $\chi\in\RPP$, 
let ${\boldsymbol B}\colon\HHH\to\HHH$ be a $\chi$-cocoercive operator,
and let $\gamma\in\left]0,2\chi\right[$. Then 
${\boldsymbol \Id}-\gamma {\boldsymbol B}$ is nonexpansive.
\end{lemma}

The main result of this section is the following theorem.

\begin{theorem}
\label{t:1}
Let $((x_{i,n})_{n\in\NN})_{1\leq i\leq m}$ be sequences generated by
Algorithm~\ref{algo:1}. Then, for every $i\in\{1,\ldots,m\}$, 
$(x_{i,n})_{n\in\NN}$ converges weakly to a point $x_i\in\HH_i$, and 
$(x_i)_{1\leq i\leq m}$ is a solution to Problem~\ref{prob:64}. 
\end{theorem}
\begin{proof}
Throughout the proof, a generic element ${\boldsymbol x}$
in the Cartesian product $\HH_1\times\cdots\times\HH_m$ will 
be expressed in terms of its components as 
${\boldsymbol x}=(x_i)_{1\leq i\leq m}$.
We shall show our algorithmic setting reduces to the situation 
described in Lemma~\ref{l:10} in
the real Hilbert space $\HHH$ obtained by endowing 
$\HH_1\times\cdots\times\HH_m$ with the scalar product 
\begin{equation}
\label{e:genna07-10}
\pscal{\cdot}{\cdot}\colon({\boldsymbol x},{\boldsymbol y})\mapsto
\sum_{i=1}^m\scal{x_i}{y_i},
\end{equation}
with associated norm 
\begin{equation}
\label{e:genna07-11}
|||\cdot|||\colon{\boldsymbol x}\mapsto
\sqrt{\sum_{i=1}^m\|x_i\|^2}.
\end{equation}
To this end, we shall show that the iterations \eqref{e:genna07-2}
can be cast in the form of \eqref{e:main2}.
First, define 
\begin{equation}
\label{e:genna07-12}
{\boldsymbol A}\colon\HHH\to 2^{\HHH}\colon{\boldsymbol x}\mapsto
\underset{i=1}{\overset{m}{\cart}}A_ix_i
\quad\text{and}\quad(\forall n\in\NN)\quad
{\boldsymbol A_n}\colon\HHH\to 2^{\HHH}\colon{\boldsymbol x}\mapsto
\underset{i=1}{\overset{m}{\cart}}A_{i,n}x_i.
\end{equation}
It follows from the maximal monotonicity of the operators 
$(A_i)_{1\leq i\leq m}$, condition~\ref{c:1} in Algorithm~\ref{algo:1}, 
\eqref{e:genna07-10}, and \eqref{e:genna07-12} that 
\begin{equation}
\label{e:genna07-18}
{\boldsymbol A}\;\text{and}\;({\boldsymbol A}_n)_{n\in\NN}
\;\text{are maximal monotone},
\end{equation}
with resolvents
\begin{equation}
\label{e:genna07-9}
J_{{\boldsymbol A}}\colon\HHH\to\HHH\colon
{\boldsymbol x}\mapsto(J_{A_i}x_i)_{1\leq i\leq m}
\quad\text{and}\quad(\forall n\in\NN)\quad
J_{{\boldsymbol A}_n}\colon\HHH\to\HHH\colon
{\boldsymbol x}\mapsto(J_{A_{i,n}}x_i)_{1\leq i\leq m},
\end{equation}
respectively. Moreover, for every $\rho\in\RPP$, we derive from 
\eqref{e:genna07-11}, \eqref{e:genna07-9}, and 
condition~\ref{c:1} in Algorithm~\ref{algo:1} that
\begin{align}
\label{e:8nov2007-20}
\sum_{n\in\NN}\sup_{|||{\boldsymbol y}|||\leq\rho}
|||J_{\gamma_n{\boldsymbol A}_n}{\boldsymbol y}
-J_{\gamma_n{\boldsymbol A}}{\boldsymbol y}|||
&=\sum_{n\in\NN}\sup_{|||{\boldsymbol y}|||\leq\rho}\sqrt{\sum_{i=1}^m
\|J_{\gamma_nA_{i,n}}y_i-J_{\gamma_nA_i}y_i\|^2}
\nonumber\\
&\leq\sum_{n\in\NN}\sup_{|||{\boldsymbol y}|||\leq\rho}\:\sum_{i=1}^m
\|J_{\gamma_nA_{i,n}}y_i-J_{\gamma_nA_i}y_i\|
\nonumber\\
&\leq\sum_{i=1}^m\sum_{n\in\NN}\:\sup_{\|y_i\|\leq\rho}
\|J_{\gamma_nA_{i,n}}y_i-J_{\gamma_nA_i}y_i\|
\nonumber\\
&<\pinf.
\end{align}
Now define
\begin{equation}
\label{e:genna07-12+}
{\boldsymbol B}\colon\HHH\to\HHH\colon{\boldsymbol x}\mapsto
(B_{i}{\boldsymbol x})_{1\leq i\leq m}
\quad\text{and}\quad(\forall n\in\NN)\quad
{\boldsymbol B_{n}}\colon\HHH\to\HHH\colon{\boldsymbol x}\mapsto
(B_{i,n}{\boldsymbol x})_{1\leq i\leq m}.
\end{equation}
Then \eqref{e:genna07-1} is equivalent to
\begin{equation}
\label{e:genna07-8}
\text{find}\quad{\boldsymbol x}\in\HHH\quad\text{such that}
\quad{\boldsymbol  0}\in{\boldsymbol A}{\boldsymbol x}
+{\boldsymbol B}{\boldsymbol x}.
\end{equation}
Moreover, in the light of \eqref{e:genna07-10}, \eqref{e:genna07-11}, and
\eqref{e:genna07-12+}, \eqref{e:genna07-21} becomes 
\begin{equation}
\label{e:genna07-7}
(\forall{\boldsymbol x}\in\HHH)
(\forall{\boldsymbol y}\in\HHH)\quad
\pscal{{\boldsymbol x}-{\boldsymbol y}}
{{\boldsymbol B}{\boldsymbol x}-{\boldsymbol B}{\boldsymbol y}}
\geq\beta|||{\boldsymbol B}{\boldsymbol x}-{\boldsymbol B}
{\boldsymbol y}|||^2. 
\end{equation}
In other words, ${\boldsymbol B}$ is 
$\beta$-cocoercive. Next, let $n\in\NN$ and set
\begin{equation}
\label{e:genna07-15}
{\boldsymbol c}_n=(a_{i,n})_{1\leq i\leq m}
\quad\text{and}\quad
{\boldsymbol d}_n=(b_{i,n})_{1\leq i\leq m}.
\end{equation}
We deduce from \eqref{e:genna07-11} and condition~\ref{c:3} 
in Algorithm~\ref{algo:1} that
\begin{equation}
\label{e:genna07-14}
\sum_{k\in\NN}|||{\boldsymbol c}_k|||\leq
\sum_{k\in\NN}\sqrt{\sum_{i=1}^m\|a_{i,k}\|^2}\leq
\sum_{i=1}^m\sum_{k\in\NN}\|a_{i,k}\|<\pinf
\end{equation}
and, likewise, that
\begin{equation}
\label{e:genna07-14+}
\sum_{k\in\NN}|||{\boldsymbol d}_k|||<\pinf. 
\end{equation}
Now set
\begin{equation}
\label{e:genna07-15+}
{\boldsymbol x}_n=(x_{i,n})_{1\leq i\leq m}
\quad\text{and}\quad
{\boldsymbol\Lambda}_n\colon\HHH\to\HHH\colon 
{\boldsymbol x}\mapsto(\lambda_{i,n}x_i)_{1\leq i\leq m}.
\end{equation}
It follows from \eqref{e:genna07-11} and condition~\ref{c:4} in 
Algorithm~\ref{algo:1} that 
\begin{equation}
\label{e:7nov2007-1}
|||{\boldsymbol\Lambda}_n|||=\max_{1\leq i\leq m}\lambda_{i,n}\leq 1
\quad\text{and}\quad
|||{\boldsymbol\Id}-{\boldsymbol\Lambda}_n|||=
1-\min_{1\leq i\leq m}\lambda_{i,n}\leq 1.
\end{equation}
Hence,
\begin{equation}
\label{e:7nov2007-2}
|||{\boldsymbol\Lambda}_n|||+
|||{\boldsymbol\Id}-{\boldsymbol\Lambda}_n|||=1+
\max_{1\leq i\leq m}(\lambda_{i,n}-\lambda_n)-
\min_{1\leq i\leq m}(\lambda_{i,n}-\lambda_n)\leq 1+\tau_n,
\end{equation}
where
\begin{equation}
\label{e:7nov2007-3}
\tau_n=2\max_{1\leq i\leq m}|\lambda_{i,n}-\lambda_n|.
\end{equation}
We observe that, by virtue of condition~\ref{c:4} in 
Algorithm~\ref{algo:1},
\begin{equation}
\label{e:8nov2007-21}
\sum_{k\in\NN}\tau_k=
2\sum_{k\in\NN}\max_{1\leq i\leq m}|\lambda_{i,k}-\lambda_k|\leq
2\sum_{i=1}^m\sum_{k\in\NN}|\lambda_{i,k}-\lambda_k|<\pinf.
\end{equation}
Moreover, in view of \eqref{e:genna07-9}, \eqref{e:genna07-12+}, 
\eqref{e:genna07-15}, and \eqref{e:genna07-15+},
the iterations \eqref{e:genna07-2} are equivalent to
\begin{equation}
\label{e:main71}
{\boldsymbol x}_{n+1}
={\boldsymbol\Lambda}_n{\boldsymbol x}_n+
({\boldsymbol\Id}-{\boldsymbol\Lambda}_n)
\big(J_{\gamma_n {\boldsymbol A}_n}
\big({\boldsymbol x}_n-\gamma_n({\boldsymbol B}_n{\boldsymbol x}_n+
{\boldsymbol d}_n)\big)+{\boldsymbol c}_n\big).
\end{equation}
Now define
\begin{equation}
\label{e:genna07-15++}
{\boldsymbol D}_n={\boldsymbol B}_n-{\boldsymbol B}.
\end{equation}
It follows from condition~\ref{c:21} in Algorithm~\ref{algo:1}, 
\eqref{e:genna07-11}, and \eqref{e:genna07-12+} that 
${\boldsymbol D}_n$ is Lipschitz continuous with constant
$\kappa_n=\sqrt{\sum_{i=1}^m\kappa_{i,n}^2}$ and that 
\begin{equation}
\label{e:8nov2007-24}
\sum_{k\in\NN}\kappa_k=\sum_{k\in\NN}\sqrt{\sum_{i=1}^m\kappa_{i,k}^2}
\leq\sum_{i=1}^m\sum_{k\in\NN}\kappa_{i,k}<\pinf.
\end{equation}
Furthermore, set
\begin{equation}
\label{e:genna07-15--}
{\boldsymbol b}_n={\boldsymbol D}_n{\boldsymbol x}_n+{\boldsymbol d}_n
\end{equation}
and let ${\boldsymbol x}$ be a solution to Problem~\ref{prob:64}. Then
\begin{align}
\label{e:6nov2007-4}
|||{\boldsymbol b}_n|||
&\leq|||{\boldsymbol D}_n{\boldsymbol x}_n|||
+|||{\boldsymbol d}_n|||\nonumber\\
&\leq|||{\boldsymbol D}_n{\boldsymbol x}_n-
{\boldsymbol D}_n{\boldsymbol x}|||+
|||{\boldsymbol D}_n{\boldsymbol x}-{\boldsymbol D}_n{\boldsymbol z}|||
+|||{\boldsymbol d}_n|||\nonumber\\
&\leq\kappa_n(|||{\boldsymbol x}_n-{\boldsymbol x}|||+
|||{\boldsymbol x}-{\boldsymbol z}|||)+|||{\boldsymbol d}_n|||,
\end{align}
where ${\boldsymbol z}$ is provided by assumption~\ref{c:22} in
Algorithm~\ref{algo:1}. We now set
\begin{equation}
\label{e:9nov2007-1}
{\boldsymbol T}_n={\boldsymbol\Id}-\gamma_n{\boldsymbol B}
\quad\text{and}\quad
{\boldsymbol e}_n=J_{\gamma_n{\boldsymbol A}_n}
\big({\boldsymbol T}_n{\boldsymbol x}\big)-{\boldsymbol x}.
\end{equation}
On the one hand, the inequality $\sup_{k\in\NN}\gamma_k\leq 2\beta$ yields
\begin{equation}
\label{e:6nov2007-7-}
|||{\boldsymbol T}_n{\boldsymbol x}|||\leq\rho,\quad\text{where}\quad 
\rho=|||{\boldsymbol x}|||+2\beta|||{\boldsymbol B}{\boldsymbol x}|||.
\end{equation}
On the other hand, Proposition~\ref{p:10octobre2007} and 
\eqref{e:genna07-9} supply
\begin{equation}
\label{e:6nov2007-1}
{\boldsymbol x}=J_{\gamma_n{\boldsymbol A}}
({\boldsymbol T}_n{\boldsymbol x}).
\end{equation}
Therefore, \eqref{e:9nov2007-1}, \eqref{e:6nov2007-7-}, and 
\eqref{e:8nov2007-20} imply that
\begin{equation}
\label{e:6nov2007-5}
\sum_{k\in\NN}|||{\boldsymbol e}_k|||
=\sum_{k\in\NN}|||J_{\gamma_k{\boldsymbol A}_k}
({\boldsymbol T}_k{\boldsymbol x})-{\boldsymbol x}|||
=\sum_{k\in\NN}|||J_{\gamma_k{\boldsymbol A}_k}
({\boldsymbol T}_k{\boldsymbol x})
-J_{\gamma_k{\boldsymbol A}}({\boldsymbol T}_k{\boldsymbol x})|||<\pinf.
\end{equation}
In addition, \eqref{e:genna07-15++}, \eqref{e:genna07-15--}, and
\eqref{e:9nov2007-1} yield
\begin{equation}
\label{e:main72}
J_{\gamma_n{\boldsymbol A}_n}
\big({\boldsymbol x}_n-\gamma_n({\boldsymbol B}_n{\boldsymbol x}_n+
{\boldsymbol d}_n)\big)-{\boldsymbol x}
=J_{\gamma_n{\boldsymbol A}_n}
\big({\boldsymbol T}_n{\boldsymbol x}_n-\gamma_n{\boldsymbol b}_n\big)
-J_{\gamma_n{\boldsymbol A}_n}
({\boldsymbol T}_n{\boldsymbol x})+{\boldsymbol e}_n.
\end{equation}
Since $J_{\gamma_n{\boldsymbol A}}$ is nonexpansive as a 
resolvent (see \cite[Proposition~3.5.3]{Aubi90} or
\cite[Proposition~2.2.iii)]{Brez73})
and ${\boldsymbol T}_n$ is nonexpansive by Lemma~\ref{l:9},
we derive from \eqref{e:main72} and \eqref{e:6nov2007-4} that 
\begin{align}
\label{e:main75}
|||J_{\gamma_n{\boldsymbol A}_n}
\big({\boldsymbol x}_n-\gamma_n({\boldsymbol B}_n{\boldsymbol x}_n+
{\boldsymbol d}_n)\big)-{\boldsymbol x}|||
&\leq|||J_{\gamma_n{\boldsymbol A}_n}
\big({\boldsymbol T}_n{\boldsymbol x}_n-\gamma_n{\boldsymbol b}_n\big)
-J_{\gamma_n{\boldsymbol A}_n}
({\boldsymbol T}_n{\boldsymbol x})|||
+|||{\boldsymbol e}_n|||\nonumber\\
&\leq|||{\boldsymbol T}_n{\boldsymbol x}_n-\gamma_n{\boldsymbol b}_n
-{\boldsymbol T}_n{\boldsymbol x}|||
+|||{\boldsymbol e}_n|||\nonumber\\
&\leq|||{\boldsymbol x}_n-{\boldsymbol x}|||+
\gamma_n|||{\boldsymbol b}_n|||
+|||{\boldsymbol e}_n|||\nonumber\\
&\leq|||{\boldsymbol x}_n-{\boldsymbol x}|||+
2\beta|||{\boldsymbol b}_n|||
+|||{\boldsymbol e}_n|||\nonumber\\
&\leq(1+2\beta\kappa_n)|||{\boldsymbol x}_n-{\boldsymbol x}|||
+2\beta\kappa_n|||{\boldsymbol x}-{\boldsymbol z}|||\nonumber\\
&\quad\;+2\beta|||{\boldsymbol d}_n|||+|||{\boldsymbol e}_n|||.
\end{align}
Thus, it results from \eqref{e:main71}, \eqref{e:main75}, 
\eqref{e:7nov2007-2}, and \eqref{e:7nov2007-1} that
\begin{align}
\label{e:8nov2007-2}
|||{\boldsymbol x}_{n+1}-{\boldsymbol x}|||
&=|||{\boldsymbol\Lambda}_n({\boldsymbol x}_n-{\boldsymbol x})+
({\boldsymbol\Id}-{\boldsymbol\Lambda}_n)
\big(J_{\gamma_n {\boldsymbol A}_n}
\big({\boldsymbol x}_n-\gamma_n({\boldsymbol B}_n{\boldsymbol x}_n+
{\boldsymbol d}_n)\big)-{\boldsymbol x}+{\boldsymbol c}_n\big)|||
\nonumber\\
&\leq|||{\boldsymbol\Lambda}_n|||\,|||{\boldsymbol x}_n-{\boldsymbol x}|||
+|||{\boldsymbol\Id}-{\boldsymbol\Lambda}_n|||\,
|||{\boldsymbol c}_n|||
\nonumber\\
&\quad\;+
|||{\boldsymbol\Id}-{\boldsymbol\Lambda}_n|||\,
|||J_{\gamma_n {\boldsymbol A}_n}
\big({\boldsymbol x}_n-\gamma_n({\boldsymbol B}_n{\boldsymbol x}_n+
{\boldsymbol d}_n)\big)-{\boldsymbol x}|||  
\nonumber\\
&\leq|||{\boldsymbol\Lambda}_n|||\,|||{\boldsymbol x}_n-{\boldsymbol x}|||
+|||{\boldsymbol\Id}-{\boldsymbol\Lambda}_n|||\,
|||{\boldsymbol c}_n||| \nonumber\\
&\quad\;+
|||{\boldsymbol\Id}-{\boldsymbol\Lambda}_n|||\,
\big((1+2\beta\kappa_n)|||{\boldsymbol x}_n-{\boldsymbol x}|||
+2\beta\kappa_n|||{\boldsymbol x}-{\boldsymbol z}|||\nonumber\\
&\quad\;+2\beta|||{\boldsymbol d}_n|||+|||{\boldsymbol e}_n|||\big)  
\nonumber\\
&\leq\big(|||{\boldsymbol\Lambda}_n|||+
|||{\boldsymbol\Id}-{\boldsymbol\Lambda}_n|||\big)
|||{\boldsymbol x}_n-{\boldsymbol x}|||+
|||{\boldsymbol\Id}-{\boldsymbol\Lambda}_n|||\,
\big(|||{\boldsymbol c}_n|||
+2\beta\kappa_n|||{\boldsymbol x}_n-{\boldsymbol x}|||\nonumber\\
&\quad\;
+2\beta\kappa_n|||{\boldsymbol x}-{\boldsymbol z}|||
+2\beta|||{\boldsymbol d}_n|||+|||{\boldsymbol e}_n|||\big)  
\nonumber\\
&\leq(1+\tau_n)|||{\boldsymbol x}_n-{\boldsymbol x}|||+
|||{\boldsymbol c}_n|||
+2\beta\kappa_n|||{\boldsymbol x}_n-{\boldsymbol x}|||\nonumber\\
&\quad\;
+2\beta\kappa_n|||{\boldsymbol x}-{\boldsymbol z}|||
+2\beta|||{\boldsymbol d}_n|||+|||{\boldsymbol e}_n|||
\nonumber\\
&\leq(1+\alpha_n)|||{\boldsymbol x}_n-{\boldsymbol x}|||+\delta_n,
\end{align}
where
\begin{equation}
\label{e:7nov2007-5}
\alpha_n=\tau_n+2\beta\kappa_n\quad\text{and}\quad
\delta_n=|||{\boldsymbol c}_n|||+
2\beta\kappa_n|||{\boldsymbol x}-{\boldsymbol z}|||
+2\beta|||{\boldsymbol d}_n|||+|||{\boldsymbol e}_n|||.
\end{equation}
In turn, it follows from \eqref{e:8nov2007-21}, \eqref{e:8nov2007-24},
\eqref{e:genna07-14}, \eqref{e:genna07-14+}, and \eqref{e:6nov2007-5} 
that $\sum_{k\in\NN}\alpha_k<\pinf$ and $\sum_{k\in\NN}\delta_k<\pinf$.
Thus, \eqref{e:8nov2007-2} and \cite[Lemma~2.2.2]{Poly87} yield
\begin{equation}
\label{e:6nov2007-7}
\sup_{k\in\NN}|||{\boldsymbol x}_k-{\boldsymbol x}|||<\pinf
\end{equation}
and, using \eqref{e:8nov2007-24} and
\eqref{e:genna07-14+}, we derive from \eqref{e:6nov2007-4} that
\begin{equation}
\label{e:8nov2007-4}
\sum_{k\in\NN}|||{\boldsymbol b}_k|||<\pinf.
\end{equation}
In view of \eqref{e:genna07-15--} and \eqref{e:9nov2007-1}, 
\eqref{e:main71} is equivalent to
\begin{align}
\label{e:main76}
{\boldsymbol x}_{n+1}
&={\boldsymbol\Lambda}_n{\boldsymbol x}_n+
({\boldsymbol\Id}-{\boldsymbol\Lambda}_n)
\big(J_{\gamma_n{\boldsymbol A}}({\boldsymbol T}_n
{\boldsymbol x}_n-\gamma_n{\boldsymbol b}_n)+{\boldsymbol h}_n\big),
\end{align}
where
\begin{equation}
\label{e:8nov2007-1}
{\boldsymbol h}_n=
J_{\gamma_n{\boldsymbol A}_n}
({\boldsymbol T}_n{\boldsymbol x}_n-\gamma_n{\boldsymbol b}_n)
-J_{\gamma_n{\boldsymbol A}}({\boldsymbol T}_n
{\boldsymbol x}_n-\gamma_n{\boldsymbol b}_n)+{\boldsymbol c}_n.
\end{equation}
Now set $\mu=\sup_{k\in\NN}|||{\boldsymbol x}_k-{\boldsymbol x}|||+\rho
+2\beta\sup_{k\in\NN}|||{\boldsymbol b}_k|||$. Then it follows from
\eqref{e:6nov2007-7}, and \eqref{e:8nov2007-4} that $\mu<\pinf$.
Moreover, we deduce from the nonexpansivity of ${\boldsymbol T}_n$ and
\eqref{e:6nov2007-7-} that
\begin{align}
|||{\boldsymbol T}_n{\boldsymbol x}_n-\gamma_n{\boldsymbol b}_n|||
&\leq
|||{\boldsymbol T}_n{\boldsymbol x}_n-{\boldsymbol T}_n{\boldsymbol x}|||+
|||{\boldsymbol T}_n{\boldsymbol x}|||+
2\beta|||{\boldsymbol b}_n|||\nonumber\\
&\leq|||{\boldsymbol x}_n-{\boldsymbol x}|||+\rho+
2\beta|||{\boldsymbol b}_n|||\nonumber\\
&\leq\mu.
\end{align}
Hence, appealing to \eqref{e:8nov2007-20} and \eqref{e:genna07-14}, 
we deduce from \eqref{e:8nov2007-1} that
\begin{equation}
\label{e:8nov2007-5}
\sum_{k\in\NN}|||{\boldsymbol h}_k|||<\pinf.
\end{equation}
Note that, upon introducing
\begin{equation}
\label{e:8nov2007-7}
{\boldsymbol a}_n={\boldsymbol h}_n
+\frac{1}{1-\lambda_n}({\boldsymbol\Lambda}_n-\lambda_n{\boldsymbol\Id})
\big({\boldsymbol x}_n-J_{\gamma_n{\boldsymbol A}}({\boldsymbol T}_n
{\boldsymbol x}_n-\gamma_n{\boldsymbol b}_n)-{\boldsymbol h}_n\big),
\end{equation}
we can rewrite \eqref{e:main76} in the equivalent form \eqref{e:main2}, 
namely
\begin{align}
\label{e:main77}
{\boldsymbol x}_{n+1}
&=\lambda_n{\boldsymbol x}_n+(1-\lambda_n)
\big(J_{\gamma_n{\boldsymbol A}}({\boldsymbol x}_n-\gamma_n({\boldsymbol B}
{\boldsymbol x}_n+{\boldsymbol b}_n))+{\boldsymbol a}_n\big).
\end{align}
Using \eqref{e:6nov2007-1} and the nonexpansivity of 
$J_{\gamma_n{\boldsymbol A}}$ and ${\boldsymbol T}_n$, we get
\begin{align}
\label{e:8nov2007-9}
|||{\boldsymbol x}_n-J_{\gamma_n{\boldsymbol A}}({\boldsymbol T}_n
{\boldsymbol x}_n-\gamma_n{\boldsymbol b}_n)-{\boldsymbol h}_n|||
&\leq
|||{\boldsymbol x}_n-{\boldsymbol x}|||+
|||J_{\gamma_n{\boldsymbol A}}({\boldsymbol T}_n{\boldsymbol x})-
J_{\gamma_n{\boldsymbol A}}({\boldsymbol T}_n
{\boldsymbol x}_n-\gamma_n{\boldsymbol b}_n)|||\nonumber\\
&\quad\;+|||{\boldsymbol h}_n|||
\nonumber\\
&\leq
2|||{\boldsymbol x}_n-{\boldsymbol x}|||+
2\beta|||{\boldsymbol b}_n|||+|||{\boldsymbol h}_n|||.
\end{align}
Therefore, we derive from \eqref{e:6nov2007-7}, \eqref{e:8nov2007-4}, 
and \eqref{e:8nov2007-5} that
\begin{equation}
\nu=\sup_{k\in\NN}|||{\boldsymbol x}_k-
J_{\gamma_k{\boldsymbol A}}({\boldsymbol T}_k{\boldsymbol x}_k
-\gamma_k{\boldsymbol b}_k)-{\boldsymbol h}_k|||<\pinf, 
\end{equation}
and hence, from \eqref{e:8nov2007-7} and the inequality
$\lambda_n\leq 1-\varepsilon$, that
\begin{align}
\label{e:8nov2007-8}
|||{\boldsymbol a}_n|||
&\leq|||{\boldsymbol h}_n|||+\frac{1}{1-\lambda_n}
|||{\boldsymbol\Lambda}_n-\lambda_n{\boldsymbol\Id}|||
\;|||{\boldsymbol x}_n-J_{\gamma_n{\boldsymbol A}}({\boldsymbol T}_n
{\boldsymbol x}_n-\gamma_n{\boldsymbol b}_n)-{\boldsymbol h}_n|||
\nonumber\\
&\leq|||{\boldsymbol h}_n|||
+\frac{\nu}{\varepsilon}\max_{1\leq i\leq m}|\lambda_{i,n}-\lambda_n|.
\end{align}
Thus, using \eqref{e:8nov2007-5} and arguing as in
\eqref{e:8nov2007-21}, we get
\begin{equation}
\label{e:8nov2007-18}
\sum_{k\in\NN}|||{\boldsymbol a}_k|||<\pinf.
\end{equation}
However, Lemma~\ref{l:10} asserts that, under \eqref{e:genna07-18}, 
\eqref{e:genna07-7}, \eqref{e:8nov2007-4}, 
\eqref{e:8nov2007-18}, and the hypotheses on $(\gamma_n)_{n\in\NN}$ and
$(\lambda_n)_{n\in\NN}$ in Algorithm~\ref{algo:1}, the sequence 
$({\boldsymbol x}_n)_{n\in\NN}$ 
generated by \eqref{e:main77} converges weakly to a solution to
\eqref{e:genna07-8}. Since \eqref{e:main77} is equivalent to
\eqref{e:genna07-2}, and \eqref{e:genna07-8} is equivalent to 
\eqref{e:genna07-1}, the proof is complete.
\end{proof}

\section{Specific scenarios}
\label{sec:2+}

Problem~\ref{prob:64} covers various scenarios, depending on the type of 
operators $(A_i)_{1\leq i\leq m}$ utilized in \eqref{e:genna07-1}. We now 
provide some specific examples which will serve as a basis for the 
concrete problems to be discussed in Sections~\ref{sec:6}--\ref{sec:7}.

\begin{example}
\label{ex:3}
Suppose that, for every $i\in\{1,\ldots,m\}$, $A_i=\partial f_i$ where
$f_i\in\Gamma_0(\HH_i)$. Then \eqref{e:genna07-1} reduces to the system 
of coupled variational inequalities
\begin{multline}
\label{e:genna07-24}
\text{find}\;x_1\in\HH_1,\ldots,x_m\in\HH_m\quad\text{such that}\\
(\forall i\in\{1,\ldots,m\})(\forall y\in\HH_i)\quad 
\scal{x_i-y}{B_{i}(x_1,\ldots,x_m)}+f_i(x_i)\leq f_i(y).
\end{multline}
A particular case of this type of problem will be investigated in 
Section~\ref{sec:5}.
\end{example}

\begin{example}
\label{ex:4}
Suppose that, for every $i\in\{1,\ldots,m\}$, $A_i$ 
is the normal cone operator to a nonempty closed convex subset 
$C_i$ of $\HH_i$, that is
\begin{equation} 
\label{e:normalcone}
A_i=N_{C_i}\colon\HH_i\to 2^{\HH_i}\colon x\mapsto
\begin{cases}
\menge{u\in\HH}{\sup_{y\in C_i}\scal{y-x}{u}\leq 0},&\text{if}\;\;
x\in C_i;\\
\emp,&\text{otherwise}.
\end{cases}
\end{equation}
Then \eqref{e:genna07-1} becomes a system of coupled 
variational inequalities of the form
\begin{multline}
\label{e:genna07-25}
\text{find}\;x_1\in C_1,\ldots,x_m\in C_m\quad\text{such that}\\
(\forall i\in\{1,\ldots,m\})(\forall y\in C_i)\quad 
\scal{x_i-y}{B_{i}(x_1,\ldots,x_m)}\leq 0.
\end{multline}
Such formulations will be investigated in Example~\ref{ex:transport} and 
Example~\ref{ex:best approx}.
\end{example}

\begin{example}
\label{ex:5}
For every $i\in\{1,\ldots,m\}$, let ${\mathcal Y}_i$ and ${\mathcal Z}_i$ 
are real Hilbert spaces, and suppose that
\begin{equation}
\label{e:10fevrier2009}
A_i\colon{\mathcal Y}_i\oplus {\mathcal Z}_i\to 2^{{\mathcal Y}_i
\oplus {\mathcal Z}_i}\colon(y,z)\mapsto
\menge{(u,v)\in{\mathcal Y}_i\oplus {\mathcal Z}_i}
{u\in\partial(-F_i(\cdot,z))(y)
\;\;\text{and}\;\;v\in\partial(F_i(y,\cdot))(z)},
\end{equation}
where $F_i\colon{\mathcal Y}_i\oplus {\mathcal Z}_i\to\RXX$ satisfies
\begin{enumerate}
\item
\label{e:saddleK}
$(\exists (y_0,z_0)\in{\mathcal Y}_i\oplus {\mathcal Z}_i)
(\forall(y,z)\in{\mathcal Y}_i\oplus {\mathcal Z}_i)$
$F_i(y_0,z)>\minf$ and $F_i(y,z_0)<\pinf$;
\item
\label{e:saddleK'}
for every $(y,z)\in{\mathcal Y}_i\oplus {\mathcal Z}_i$, $-F_i(\cdot,z)$ 
and $F_i(y,\cdot)$ are lower semicontinuous and convex. 
\end{enumerate}
Then, for every $i\in\{1,\ldots,m\}$, $A_i$ is a maximal monotone
operator acting on $\HH_i={\mathcal Y}_i\oplus {\mathcal Z}_i$
\cite{Rock70} and, upon setting $B_i=(B_{i1},B_{i2})$, where 
$B_{i1}\colon\HH_i\to{\mathcal Y}_i$ and
$B_{i2}\colon\HH_i\to{\mathcal Z}_i$,
\eqref{e:genna07-1} reduces to the system of coupled saddle-point problems
\begin{multline}
\label{e:genna07-25o}
\text{find}\;x_1=(y_1,z_1)\in\HH_1,\:\ldots,\:x_m=(y_m,z_m)\in\HH_m 
\;\;\text{such that}\;\;(\forall i\in\{1,\ldots,m\})\\[3mm]
\begin{cases}
\sup\limits_{y\in {\mathcal Y}_i}F_i(y,z_i)-
\scal{y}{B_{i1}(x_1,\ldots,x_m)}_{{\mathcal Y}_i}
=F_i(y_i,z_i)-\scal{y_i}{B_{i1}
(x_1,\ldots,x_m)}_{{\mathcal Y}_i}\\
\inf\limits_{z\in {\mathcal Z}_i}F_i(y_i,z)+
\scal{z}{B_{i2}(x_1,\ldots,x_m)}_{{\mathcal Z}_i}
=F_i(y_i,z_i)+\scal{z_i}{B_{i2}
(x_1,\ldots,x_m)}_{{\mathcal Z}_i}.
\end{cases}
\end{multline}
Such formulations will arise in Example~\ref{ex:saddle}.
\end{example}

\begin{example}
\label{ex:brezis}
Let us recall some standard notation \cite{Brez73,Zeid90}. 
Fix $T\in\RPP$ and $p\in\left[1,\pinf\right[$.
Then ${\EuScript D}(]0,T[)$ is the set of infinitely differentiable 
functions from $]0,T[$ to $\RR$ with compact support in $]0,T[$. Given a 
real Hilbert space ${\mathsf H}$, ${\EuScript C}([0,T];{\mathsf H})$ 
is the space of continuous functions from $[0,T]$ to ${\mathsf H}$ and
$L^p([0,T];{\mathsf H})$ is the space of 
classes of equivalences of Borel measurable functions 
$x\colon[0,T]\to{\mathsf H}$ such that
$\int_0^T\|x(t)\|_{\mathsf H}^pdt<\pinf$.
$L^2([0,T];{\mathsf H})$ is a Hilbert space with scalar product
$(x,y)\mapsto\int_0^T\scal{x(t)}{y(t)}_{\mathsf H}dt$.
Now take $x$ and $y$ in $L^1([0,T];{\mathsf H})$. Then $y$ is the
weak derivative of $x$ if 
$\int_0^T\phi(t)y(t)dt=-\int_0^T(d\phi(t)/dt)x(t)dt$ for every 
$\phi\in{\EuScript D}(]0,T[)$, in which case we use the notation 
$y=x'$. Moreover, 
\begin{equation}
W^{1,2}([0,T];{\mathsf H})=
\menge{x\in L^2([0,T];{\mathsf H})}{x'\in L^2([0,T];{\mathsf H})}.
\end{equation}
Now, for every $i\in\{1,\ldots,m\}$, let ${\mathsf H}_i$ be a 
real Hilbert space, let 
${\mathsf A_i}\colon {\mathsf H}_i\to 2^{{\mathsf H}_i}$ be 
maximal monotone, let ${\mathsf B}_i\colon{\mathsf H}_1\times\cdots\times
{\mathsf H}_m\to{\mathsf H}_i$, and set $\HH_i=L^2([0,T];{\mathsf H}_i)$.
Then, under standard assumptions, the operator 
\begin{equation}
\label{e:08-08-08'}
A_i\colon\HH_i\to 2^{{\HH}_i}
\colon x\mapsto
\begin{cases}
x'+{\mathsf A}_ix,
&\text{if}\;\;x\in W^{1,2}([0,T];{\mathsf H}_i)\;\;\text{and}\;\;
x(0)=x(T);\\
\emp,&\text{otherwise}
\end{cases}
\end{equation}
is maximal monotone (see \cite{Mich96}, \cite[Section~3.6]{Brez73}, 
\cite{Simo08}). In this context, with a suitable construction of the 
operators $(B_i)_{1\leq i\leq m}$ in terms of 
$({\mathsf B}_i)_{1\leq i\leq m}$, \eqref{e:genna07-1} assumes the 
form of the system of coupled evolution inclusions
\begin{multline}
\label{e:genna07-38}
\text{find}\;x_1\in W^{1,2}([0,T];{\mathsf H}_1),\ldots,
x_m\in W^{1,2}([0,T];{\mathsf H}_m)\quad\text{such that}\\
(\forall i\in\{1,\ldots,m\})\quad 
\begin{cases}
0\in\displaystyle{x_i'(t)}+
{\mathsf A}_i (x_i(t))+{\mathsf B}_i(x_1 (t),\ldots,x_m (t))\;\;
\text{a.e. on}\;\left]0,T\right[\\
x_i(0)=x_i(T).
\end{cases}
\end{multline}
This type of problem will be revisited in Section~\ref{sec:6}.
\end{example}

In Algorithm~\ref{algo:1}, maximal monotone approximations 
$(A_{i,n})_{1\leq i\leq m}$ to the original operators 
$(A_i)_{1\leq i\leq m}$ can be used at iteration $n$, as long as 
\eqref{e:lc1} is satisfied.
In order to illustrate this condition, we need a couple of 
definitions and some technical facts.

\begin{definition}
\label{d:h}
Let $A$ and $B$ be set-valued operators from $\HH$ to $2^{\HH}$ and let 
$\varrho\in\RPP$ be such that $E_\varrho\cap(\gr A\cup\gr B)\neq\emp$, 
where 
$E_\varrho=\menge{(x,y)\in\HH\times\HH}{\max\{\|x\|,\|y\|\}\leq\varrho}$. 
The $\varrho$-Hausdorff distance between $A$ and $B$ is 
\cite[Section~1.1]{Atto93}
\begin{equation}
\label{e:rho-haus}
\haus_\varrho(A,B)
=\max\bigg\{\sup_{z\in E_\varrho\cap\gr B}d_{\gr A}(z),
\sup_{z\in E_\varrho\cap\gr A}d_{\gr B}(z)\bigg\}.
\end{equation}
Moreover, the Yosida approximation of $A$ of index $\gamma\in\RPP$ is 
$\moyo{A}{\gamma}=(\Id-J_{\gamma A})/\gamma$ \cite{Aubi90,Brez73}.
\end{definition}

\begin{lemma}
\label{l:18}
Let $A\colon\HH\to 2^{\HH}$ be maximal monotone, let $x\in\HH$, 
let $\gamma\in\RPP$, and let $\mu\in\RPP$. Then the following hold.
\begin{enumerate}
\item
\label{l:18i}
$J_{\mu A}x=J_{\gamma A}\big(x+(1-\gamma/\mu)(J_{\mu A}x-x)\big)$.
\item
\label{l:18i+}
$\gamma\leq\mu$ $\Rightarrow$ $\|J_{\gamma A}x-x\|\leq 2\|J_{\mu A}x-x\|$.
\item
\label{l:18ii}
$J_{\gamma(\moyo{A}{\mu})}x=
x+\gamma\big(J_{(\gamma+\mu)A}x-x\big)/(\gamma+\mu)$.
\item
\label{l:18iii}
$\big\|J_{\gamma (\moyo{A}{\mu})}x-J_{\gamma A}x\big\|
\leq 2\mu\|J_{\gamma A}x-x\|/(\gamma+\mu)$.
\end{enumerate}
\end{lemma}
\begin{proof}
\ref{l:18i}: See \cite[Section~II.4]{Brez73}.

\ref{l:18i+}: 
Set $\lambda=\gamma/\mu$ and observe that $\lambda\in\left]0,1\right]$.
It follows from the nonexpansivity of $J_{\gamma A}$ 
\cite[Proposition~3.5.3]{Aubi90} and \ref{l:18i} that 
$\|J_{\gamma A}x-x\|\leq\|J_{\gamma A}x-J_{\mu A}x\|+
\|J_{\mu A}x-x\|=\|J_{\gamma A}x-J_{\gamma A}\left(\lambda x+
\left(1-\lambda\right)J_{\mu A}x\right)\|+
\|J_{\mu A}x-x\|\leq\|x-\lambda x-
\left(1-\lambda\right)J_{\mu A}x\|+
\|J_{\mu A}x-x\|\leq 2\|J_{\mu A}x-x\|$.

\ref{l:18ii}: This identity follows at once from the semigroup property
$\moyo{A}{\gamma+\mu}=\moyo{(\!\moyo{A}{\mu})}{\gamma\;\,}$,
which can be found in \cite[Proposition~2.6(ii)]{Brez73}.

\ref{l:18iii}:
It follows from \ref{l:18ii} that
\begin{align}
\label{e:develop}
\big\|J_{\gamma (\moyo{A}{\mu})}x-J_{\gamma A}x\big\|
&=\Big\|x+\frac{\gamma}{\gamma+\mu}
\Big(J_{(\gamma+\mu)A}x-x\Big)-J_{\gamma A}x\Big\|
\nonumber\\
&=\Big\|\frac{\gamma}{\gamma+\mu}
\Big(J_{(\gamma+\mu)A}x-x-(J_{\gamma A}x-x)\Big)
-\frac{\mu}{\gamma+\mu}(J_{\gamma A}x-x)\Big\|\nonumber\\
&\leq \frac{\gamma}{\gamma+\mu}\|J_{(\gamma+\mu)A}x-J_{\gamma A}x\|
+\frac{\mu}{\gamma+\mu}\|J_{\gamma A}x-x\|.
\end{align}
On the other hand, it follows from \ref{l:18i} and the nonexpansivity
of $J_{(\gamma+\mu)A}$ that
\begin{align}
\|J_{(\gamma+\mu)A}x-J_{\gamma A}x\|
&=\Big\|J_{(\gamma+\mu)A}x-J_{(\gamma+\mu) A}
\Big(x+\Big(1-\frac{\gamma+\mu}{\gamma}\Big)(J_{\gamma A}x-x)\Big)
\Big\|\nonumber\\
&\leq\frac\mu\gamma\|J_{\gamma A}x-x\|
\end{align}
which, combined with \eqref{e:develop}, yields the announced inequality.
\end{proof}

\begin{proposition}
\label{p:2008-07-09}
Let $i\in\{1,\ldots,m\}$ and let $(\gamma_n)_{n\in\NN}$ be as in 
Algorithm~\ref{algo:1}. Then condition \eqref{e:lc1} holds if one of the
following is satisfied for every $n\in\NN$.
\begin{enumerate}
\item
\label{p:2008-07-09i}
$A_{i,n}=(\gamma_{i,n}/\gamma_n)A_{i}$, where
$(\gamma_{i,n})_{n\in\NN}$ lies in $\left]0,2\beta\right[$ 
and satisfies $\sum_{n\in\NN}|\gamma_{i,n}-\gamma_n|<\pinf$. 
\item
\label{p:2008-07-09ii}
$A_{i,n}=\moyo{A_i}{\mu_{i,n}}$, 
where $(\mu_{i,n})_{n\in\NN}$ lies in $\RPP$ and satisfies
$\sum_{n\in\NN}\mu_{i,n}<\pinf$.
\item
\label{p:2008-07-09iii}
$\gamma_n=\gamma\in\left[\varepsilon,2\beta-\varepsilon\right]$, and
\begin{equation}
\label{e:9juillet2008}
\big(\forall\varrho\in\big[\,\|J_{\gamma A_i}0
\|\max\{1,1/\gamma\},\pinf\big[\,\big)
\quad\sum_{n\in\NN}\haus_{\varrho}(A_i,A_{i,n})<\pinf.
\end{equation}
\end{enumerate}
\end{proposition}
\begin{proof}
Let $\rho\in\RPP$. Since $\sup_{n\in\NN}\gamma_n\leq 2\beta$,
we derive from Lemma~\ref{l:18}\ref{l:18i+} and the 
nonexpansivity of $\Id-J_{2\beta A_i}=J_{(2\beta A_i)^{-1}}$ that
\begin{align}
(\forall n\in\NN)(\forall y\in \HH_i)\quad\|J_{\gamma_n A_i}y-y\|
&\leq 2\|J_{2\beta A_i}y-y\|\nonumber\\
&\leq 2\|(\Id-J_{2\beta A_i})y-(\Id-J_{2\beta A_i})0\|+
2\|J_{2\beta A_i}0\|
\nonumber\\
&\leq2\|y\|+2\|J_{2\beta A_i}0\|.
\label{e:19juin2008}
\end{align}
In addition, set $\mu=2\rho+2\|J_{2\beta A_i}0\|$. We now prove
assertions \ref{p:2008-07-09i}--\ref{p:2008-07-09iii}.

\ref{p:2008-07-09i}:
It follows from Lemma~\ref{l:18}\ref{l:18i} and 
the nonexpansivity of $J_{\gamma_{i,n} A_i}$ that
\begin{align}
\label{e:J-J2}
(\forall n\in\NN)(\forall y\in\HH_i)\quad
\|J_{\gamma_{i,n}A_i}y-J_{\gamma_n A_i}y\|
&=\big\|J_{\gamma_{i,n}A_i}y-J_{\gamma_{i,n}A_i}
\big(y+(1-\gamma_{i,n}/\gamma_n)(J_{\gamma_nA_i}y-y)\big)\big\|\nonumber\\
&\leq|1-\gamma_{i,n}/\gamma_n|\,\|J_{\gamma_nA_i}y-y\|.
\end{align}
Hence, in view of \eqref{e:J-J2}, \eqref{e:19juin2008}, and the 
inequality $\inf_{n\in\NN}\gamma_{n}\geq\varepsilon$ we have
\begin{equation}
\sum_{n\in\NN}\:\sup_{\|y\|\leq\rho}
\|J_{\gamma_{i,n}A_i}y-J_{\gamma_n A_i}y\|
\leq\mu\sum_{n\in\NN}|1-\gamma_{i,n}/\gamma_n|
\leq\frac{\mu}{\varepsilon}\sum_{n\in\NN}|\gamma_n-\gamma_{i,n}|<\pinf,
\end{equation}
which yields \eqref{e:lc1}.

\ref{p:2008-07-09ii}:
For every $y\in\HH_i$ such that $\|y\|\leq\rho$ and every
$n\in\NN$, Lemma~\ref{l:18}\ref{l:18iii} and \eqref{e:19juin2008} yield
\begin{equation}
\label{e:ineg}
\big\|J_{\gamma_n (\moyo{A_i}{\mu_{i,n}})}y-J_{\gamma_n A_i}y\big\|
\leq\frac{2\mu_{i,n}}{\gamma_n+\mu_{i,n}}\|J_{\gamma_n A_i}y-y\|
\leq\frac{2\mu_{i,n}}{\varepsilon}\mu.
\end{equation}
Consequently, \eqref{e:lc1} holds.

\ref{p:2008-07-09iii}:
Set $\varrho=\max\{\rho+\|J_{\gamma A_i}0\|,
(\rho+\|J_{\gamma A_i}0\|)/\gamma\}$ and let 
$E_\varrho\subset\HH_i\times\HH_i$ be as in Definition~\ref{d:h}.
It follows from \cite[Proposition~1.2]{Atto93} that 
$E_\varrho\cap\gr A_i\neq\emp$ and that
\begin{equation}
(\forall n\in\NN)\quad
\sup_{\|y\|\leq \rho}\|J_{\gamma A_{i,n}}y-J_{\gamma A_i}y\|
\leq (2+\gamma) \haus_{\varrho}(A_{i,n},A_i).
\end{equation}
Since, in view of \eqref{e:9juillet2008},
$\sum_{n\in\NN}\haus_{\varrho}(A_{i,n},A_i)<\pinf$, we conclude that
\eqref{e:lc1} holds. 
\end{proof}

\section{Coupling schemes}
\label{sec:3}

The coupling between the $m$ inclusions in 
Problem~\ref{prob:64} is determined by the operators 
$(B_i)_{1\leq i\leq m}$, which must satisfy \eqref{e:genna07-21}.
In this section, we describe various situations in
which this property holds. In each case, the value of $\beta$ in
\eqref{e:genna07-21} will be specified, as it is explicitly required in 
Algorithm~\ref{algo:1}. In this connection, the notion of cocoercivity 
(see Definition~\ref{d:coco}) will play an important role.
Examples of cocoercive operators include firmly nonexpansive operators 
(e.g., resolvents of maximal monotone operators, proximity 
operators, and projection operators onto nonempty closed convex sets).
In addition, the Yosida approximation of a maximal monotone operator 
of index 
$\chi$ is $\chi$-cocoercive \cite{Atto84} (further examples of 
cocoercive operators can be found in \cite{Zhud96}). It is clear from 
\eqref{e:cocoercive} that if $T$ is $\chi$-cocoercive, then it is 
$\chi^{-1}$-Lipschitz continuous. The next lemma, which
provides a converse implication, supplies us with another important 
instance of cocoercive operator (see also \cite{Dunn76}).

\begin{lemma}{\rm\cite[Corollaire~10]{Bail77}}
\label{l:BH}
Let $\varphi\colon\HH\to\RR$ be a differentiable convex function and let
$\tau\in\RPP$. Suppose that $\nabla\varphi$ is $\tau$-Lipschitz
continuous. Then $\nabla\varphi$ is $\tau^{-1}$-cocoercive.
\end{lemma}

\begin{lemma}
\label{l:BH2}
Let $L\in\BL(\HH)$ be a nonzero self-adjoint operator such that 
$(\forall x\in\HH)$ $\scal{Lx}{x}\geq 0$. Then $L$ is 
$\|L\|^{-1}$-cocoercive.
\end{lemma}
\begin{proof}
Set $\varphi\colon x\mapsto\scal{Lx}{x}/2$. 
Then $\varphi$ is convex and differentiable, and its gradient
$\nabla\varphi\colon x\mapsto Lx$ is $\|L\|$-Lipschitz continuous.
Hence, the assertion follows from Lemma~\ref{l:BH}.
\end{proof}

\subsection{Linear coupling}
\label{sec:31}
We examine the case in which the operators $(B_i)_{1\leq i\leq m}$
are linear, which reduces \eqref{e:genna07-21} to
\begin{equation}
\label{e:17octobre2007}
(\forall(x_1,\ldots,x_m)\in\HH_1\times\cdots\times\HH_m)\quad
\sum_{i=1}^m\scal{B_{i}(x_1,\ldots,x_m)}{x_i}\geq
\beta\sum_{i=1}^m\big\|B_{i}(x_1,\ldots,x_m)\big\|^2.
\end{equation}
We assume that, for every $i$ and $j$ in 
$\{1,\ldots,m\}$, there exists $M_{ij}\in\BL(\HH_j,\HH_i)$ such that
\begin{equation}
\label{e:sud3}
(\forall i\in\{1,\ldots,m\})\quad
B_i\colon\HH_1\times\cdots\times\HH_m\to\HH_i\colon
(x_j)_{1\leq j\leq m}\mapsto\sum_{j=1}^mM_{ij}x_j.
\end{equation}
Thus, \eqref{e:17octobre2007} is equivalent to 
\begin{equation}
\label{e:18octobre2007}
(\forall(x_1,\ldots,x_m)\in\HH_1\times\cdots\times\HH_m)\quad
\sum_{i=1}^m\sum_{j=1}^m\scal{M_{ij}x_j}{x_i}\geq
\beta\sum_{i=1}^m\bigg\|\sum_{j=1}^mM_{ij}x_j\bigg\|^2.
\end{equation}
Our objective is to determine tight values of $\beta$ for which this
inequality holds in various scenarios. As in the proof of 
Theorem~\ref{t:1}, it will be convenient to let $\HHH$ be the direct 
Hilbert sum of the spaces $(\HH_i)_{1\leq i\leq m}$ with the notation
\eqref{e:genna07-10} and \eqref{e:genna07-11}, and to set 
\begin{equation}
\label{e:guadeloupe2007-1}
{\boldsymbol B}\colon\HHH\to\HHH\colon{\boldsymbol x}\mapsto
(B_{i}{\boldsymbol x})_{1\leq i\leq m}=
\bigg(\sum_{j=1}^mM_{ij}x_j\bigg)_{1\leq i\leq m}.
\end{equation}

\begin{proposition}
\label{p:12}
Suppose that the following hold.
\begin{enumerate}
\item
\label{i:sud4}
$(\exi (i,j)\in\{1,\ldots,m\}^2)$ $M_{ij}\neq 0$.
\item
\label{i:sud1}
$(\forall(i,j)\in\{1,\ldots,m\}^2)$ $M_{ji}=M^*_{ij}$.
\item
\label{i:sud2}
$(\forall(x_1,\ldots,x_m)\in\HH_1\times\cdots\times\HH_m)$
$\sum_{i=1}^m\sum_{j=1}^m\scal{M_{ij}x_j}{x_i}\geq 0$.
\end{enumerate}
Then \eqref{e:18octobre2007} is satisfied with 
$\beta=1/|||{\boldsymbol B}|||$ and, a fortiori, with 
\begin{equation}
\label{e:genna07-49+++}
\beta=\frac{1}{\sqrt{\sum_{i=1}^m\sum_{j=1}^m\|M_{ij}\|^2}}.
\end{equation}
\end{proposition}
\begin{proof}
It follows from \ref{i:sud4} that ${\boldsymbol B}\neq{\boldsymbol 0}$ 
and from \ref{i:sud1} that ${\boldsymbol B}^*={\boldsymbol B}$. 
In addition, \eqref{e:genna07-10} and \ref{i:sud2} imply that 
$(\forall{\boldsymbol x}\in\HHH)$ 
$\pscal{{\boldsymbol B}{\boldsymbol x}}{{\boldsymbol x}}\geq 0$. 
Hence, we derive from Lemma~\ref{l:BH2} that ${\boldsymbol B}$
is $|||{\boldsymbol B}|||^{-1}$-cocoercive which, in view of 
\eqref{e:guadeloupe2007-1}, \eqref{e:genna07-10}, 
and \eqref{e:genna07-11}, can be expressed as 
\begin{equation}
\label{e:guadeloupe2007-2}
(\forall(x_1,\ldots,x_m)\in\HH_1\times\cdots\times\HH_m)\quad
\sum_{i=1}^m\sum_{j=1}^m\scal{M_{ij}x_j}{x_i}\geq
\frac{1}{|||{\boldsymbol B}|||}
\sum_{i=1}^m\bigg\|\sum_{j=1}^mM_{ij}x_j\bigg\|^2.
\end{equation}
Hence, \eqref{e:18octobre2007} holds with 
$\beta=1/|||{\boldsymbol B}|||$.
Now take ${\boldsymbol x}\in\HHH$ such that 
$|||{\boldsymbol x}|||\leq 1$. Then, \eqref{e:guadeloupe2007-1} and 
Cauchy-Schwarz yield 
\begin{align}
\label{e:guadeloupe2007-3}
|||{\boldsymbol B}{\boldsymbol x}|||^2
&=\sum_{i=1}^m\bigg\|\sum_{j=1}^mM_{ij}x_j\bigg\|^2\nonumber\\
&\leq\sum_{i=1}^m\bigg(\sum_{j=1}^m\|M_{ij}\|\,\|x_j\|\bigg)^2
\nonumber\\
&\leq\sum_{i=1}^m\bigg(\sum_{j=1}^m\|M_{ij}\|^2\bigg)
\bigg(\sum_{j=1}^m\|x_j\|^2\bigg)\nonumber\\
&\leq\sum_{i=1}^m\sum_{j=1}^m\|M_{ij}\|^2.
\end{align}
Thus, $|||{\boldsymbol B}|||^2\leq\sum_{i=1}^m\sum_{j=1}^m\|M_{ij}\|^2$
and it follows from \eqref{e:guadeloupe2007-2} that 
\eqref{e:18octobre2007} holds with \eqref{e:genna07-49+++}.
\end{proof}

In practice, one is interested in obtaining the tightest bound in 
\eqref{e:18octobre2007}. If $|||{\boldsymbol B}|||$ is known,
one will use $\beta=1/|||{\boldsymbol B}|||$ in 
Algorithm~\ref{algo:1}. This is for instance the case
in the next proposition. In many situations, however, 
$|||{\boldsymbol B}|||$ will be hard to compute and one can use the 
value supplied by \eqref{e:genna07-49+++}, which requires 
only knowledge of the norms of the individual operators 
$(M_{ij})_{1\leq i,j\leq m}$. 

\begin{proposition}
\label{p:13}
Let $\Xi=[\xi_{ij}]$ be a nonzero real $m\times m$ positive semidefinite 
symmetric matrix with largest eigenvalue $\lambda_{\operatorname{max}}$. 
Set
\begin{equation}
\label{e:sud6}
(\forall i\in\{1,\ldots,m\})\quad\HH_i=\HH\quad\text{and}\quad
B_i\colon\HH^m\to\HH\colon
(x_j)_{1\leq j\leq m}\mapsto\sum_{j=1}^m\xi_{ij}x_j.
\end{equation}
Then \eqref{e:18octobre2007} holds with 
$\beta=1/\lambda_{\operatorname{max}}$.
\end{proposition}
\begin{proof}
Let $\Lambda$ be the diagonal matrix the diagonal entries of which are
the eigenvalues $(\lambda_i)_{1\leq i\leq m}$ of $\Xi$. There exists
an $m\times m$ orthogonal matrix $\Pi=[\pi_{ij}]$ such that 
$\Xi=\Pi\Lambda\Pi^t$. Now 
set ${\boldsymbol D}\colon\HHH\to\HHH\colon{\boldsymbol x}\mapsto
\big(\lambda_ix_i\big)_{1\leq i\leq m}$ and
${\boldsymbol U}\colon\HHH\to\HHH\colon{\boldsymbol x}\mapsto
\big(\sum_{j=1}^m\pi_{ij}x_j\big)_{1\leq i\leq m}$. Then ${\boldsymbol U}$ 
is unitary and
$|||{\boldsymbol B}|||^2=|||{\boldsymbol U}{\boldsymbol D}
{\boldsymbol U}^*|||^2
=|||{\boldsymbol D}|||^2=\sup_{|||{\boldsymbol x}|||\leq 1}
\sum_{i=1}^m\lambda_i^2\|x_i\|^2=\lambda_{\max}^2$. Hence, the assertion
follows from Proposition~\ref{p:12}.
\end{proof}

As shown next, equality can be achieved in \eqref{e:17octobre2007}.
\begin{example}
Set
\begin{equation}
\label{e:22octobre2007}
(\forall i\in\{1,\ldots,m\})\quad\HH_i=\HH\quad\text{and}\quad
B_i\colon\HH^m\to\HH\colon
(x_j)_{1\leq j\leq m}\mapsto x_i-\frac1m\sum_{j=1}^mx_j.
\end{equation}
Then equality is achieved in \eqref{e:17octobre2007} with $\beta=1$. 
\end{example}
\begin{proof}
Let $(x_1,\ldots,x_m)\in\HH^m$. Then 
\begin{align}
\sum_{i=1}^m\Scal{x_i-\frac1m\sum_{j=1}^m x_j}{x_i}
&=\frac1m\bigg\|\sum_{j=1}^m x_j\bigg\|^2+
\sum_{i=1}^m\Scal{x_i-\frac2m\sum_{j=1}^m x_j}{x_i}\nonumber\\
&=\sum_{i=1}^m\bigg\|x_i-\frac1m\sum_{j=1}^m x_j\bigg\|^2,
\end{align}
which provides the announced identity. 
\end{proof}

Our last example concerns a specific structure of the operators
$(M_{ij})_{1\leq i,j\leq m}$.

\begin{proposition}
\label{p:14}
For every $k\in\{1,\ldots,p\}$, let $\GG_k$ be a real Hilbert space
and, for every $i\in\{1,\ldots,m\}$, let 
$L_{ki}\in\BL(\HH_i,\GG_k)$.
Assume that $\min_{1\leq k\leq p}\sum_{i=1}^m\|L_{ki}\|^2>0$ and set 
\begin{equation}
\label{e:4octobre2007}
(\forall(i,j)\in\{1,\ldots,m\}^2)\quad
M_{ij}=\sum_{k=1}^pL_{ki}^*L_{kj}
\end{equation}
in \eqref{e:sud3}. Then \eqref{e:18octobre2007} holds with 
\begin{equation}
\label{e:genna07-49}
\beta=\frac{1}{\sum_{k=1}^p\sum_{i=1}^m\|L_{ki}\|^2}.
\end{equation}
\end{proposition}
\begin{proof}
For every $i$ and $j$ in $\{1,\ldots,m\}$, \eqref{e:4octobre2007} and
Cauchy-Schwarz yield
\begin{equation}
\label{e:22nov2007-1}
\|M_{ij}\|^2=\bigg\|\sum_{k=1}^pL_{ki}^*L_{kj}\bigg\|^2
\leq\bigg(\sum_{k=1}^p\|L_{ki}\|\,\|L_{kj}\|\bigg)^2
\leq\bigg(\sum_{k=1}^p\|L_{ki}\|^2\bigg)
\bigg(\sum_{k=1}^p\|L_{kj}\|^2\bigg).
\end{equation}
Consequently, 
\begin{equation}
\label{e:22nov2007-2}
\sum_{i=1}^m\sum_{j=1}^m\|M_{ij}\|^2\leq
\bigg(\sum_{k=1}^p\sum_{i=1}^m\|L_{ki}\|^2\bigg)
\bigg(\sum_{k=1}^p\sum_{j=1}^m\|L_{kj}\|^2\bigg)=
\bigg(\sum_{k=1}^p\sum_{i=1}^m\|L_{ki}\|^2\bigg)^2.
\end{equation}
On the other hand, it follows from \eqref{e:4octobre2007} that
conditions~\ref{i:sud4}--\ref{i:sud2} in 
Proposition~\ref{p:12} are satisfied. Therefore, we derive from
Proposition~\ref{p:12} that \eqref{e:18octobre2007} holds with 
$\beta$ as defined in \eqref{e:genna07-49}. 
\end{proof}

\begin{remark}
\label{r:82}
For every $i\in\{1,\ldots,m\}$ and $n\in\NN$, suppose that 
$B_{i,n}\in\BL(\HHH,\HH_i)$ in Algorithm~\ref{algo:1}, say
\begin{equation}
\label{e:sud54}
B_{i,n}\colon\HHH\to\HH_i\colon (x_j)_{1\leq j\leq m}\mapsto
\sum_{j=1}^mM_{ij,n}x_j,\;\;\text{where}\;\;
(\forall j\in\{1,\ldots,m\})\;M_{ij,n}\in\BL(\HH_j,\HH_i).
\end{equation}
Then assumption~\ref{c:22} in Algorithm~\ref{algo:1} is satisfied with 
${\boldsymbol z}={\boldsymbol 0}$. In addition, suppose that 
\begin{equation}
\label{e:anr}
\max_{1\leq i\leq m}
\sum_{n\in\NN}\sqrt{\sum_{j=1}^m\|M_{ij,n}-M_{ij}\|^2}<\pinf. 
\end{equation}
Then assumption~\ref{c:21} in Algorithm~\ref{algo:1} is satisfied. Indeed, 
let ${\boldsymbol x}\in\HHH$, $i\in\{1,\ldots,m\}$, and $n\in\NN$, and
set $\kappa_{i,n}=\sqrt{\sum_{j=1}^m\|M_{ij,n}-M_{ij}\|^2}$. Then, by 
Cauchy-Schwarz,
\begin{equation}
\|(B_{i,n}-B_i){\boldsymbol x}\|=
\bigg\|\sum_{j=1}^m(M_{ij,n}-M_{ij})x_j\bigg\|
\leq\sum_{j=1}^m\|M_{ij,n}-M_{ij}\|\,\|x_j\|
\leq\kappa_{i,n}|||{\boldsymbol x}|||,
\end{equation}
where \eqref{e:anr} yields $\sum_{n\in\NN}\kappa_{i,n}<\pinf$, 
as desired.
\end{remark}

\subsection{Nonlinear coupling}
\label{sec:32}
In this section we turn our attention to the determination of 
the parameter $\beta$ in \eqref{e:genna07-21} when the operators
$(B_i)_{1\leq i\leq m}$ are nonlinear. Our first model is a nonlinear
version of Proposition~\ref{p:14}.

\begin{proposition}
\label{p:coco}
For every $k\in\{1,\ldots,p\}$, let $\GG_k$ be a real Hilbert space, 
let $\beta_k\in\RPP$, let $T_k\colon\GG_k\to\GG_k$ be 
$\beta_k$-cocoercive, and, for every $i\in\{1,\ldots,m\}$, let 
$L_{ki}\in\BL(\HH_i,\GG_k)$. Assume that 
$\min_{1\leq k\leq p}\sum_{i=1}^m\|L_{ki}\|^2>0$ and set 
\begin{equation}
\label{e:4oct2007}
(\forall i\in\{1,\ldots,m\})\quad
B_i\colon\HH_1\times\cdots\times\HH_m\to\HH_i\colon
(x_j)_{1\leq j\leq m}\mapsto\sum_{k=1}^pL_{ki}^*T_k
\bigg(\sum_{j=1}^mL_{kj}x_j\bigg).
\end{equation}
Then \eqref{e:genna07-21} holds with 
\begin{equation}
\label{e:genna07-47}
\beta=\frac1p\:\min_{1\leq k\leq p}\:
\frac{\beta_k}{\sum_{i=1}^m\|L_{ki}\|^2}.
\end{equation}
\end{proposition}
\begin{proof}
For every $i\in\{1,\ldots,m\}$, let $x_i$ and $y_i$ be points in $\HH_i$.
It follows from \eqref{e:4oct2007}, \eqref{e:genna07-47}, and the 
convexity of 
$\|\cdot\|^2$ that
\begin{multline} 
\label{e:genna07-45}
\sum_{i=1}^m\scal{B_i(x_1,\ldots,x_m)-B_i(y_1,\ldots,y_m)}{x_i-y_i}\\
\begin{aligned}[b]
&=\sum_{i=1}^m\sum_{k=1}^p\Scal{L_{ki}^*\bigg(T_k\bigg(\sum_{j=1}^m
L_{kj}x_j\bigg)-T_k\bigg(\sum_{j=1}^mL_{kj}y_j\bigg)\bigg)}
{x_i-y_i}\\
&=\sum_{i=1}^m\sum_{k=1}^p\Scal{T_k\bigg(\sum_{j=1}^m
L_{kj}x_j\bigg)-T_k\bigg(\sum_{j=1}^mL_{kj}y_j\bigg)}
{L_{ki}({x_i}-{y_i})}\\
&=\sum_{k=1}^p\Scal{T_k\bigg(\sum_{j=1}^m
L_{kj}x_j\bigg)-T_k\bigg(\sum_{j=1}^mL_{kj}y_j\bigg)}
{\sum_{i=1}^mL_{ki}{x_i}-\sum_{i=1}^mL_{ki}{y_i}}\\
&\geq\sum_{k=1}^p\beta_k\bigg\|T_k\bigg(\sum_{j=1}^m
L_{kj}x_j\bigg)-T_k\bigg(\sum_{j=1}^mL_{kj}y_j\bigg)\bigg\|^2\\
&=\sum_{k=1}^p\frac{\beta_k}{\sum_{i=1}^m\|L_{ki}\|^2}
\sum_{i=1}^m\|L_{ki}\|^2\,\bigg\|T_k\bigg(\sum_{j=1}^m
L_{kj}x_j\bigg)-T_k\bigg(\sum_{j=1}^mL_{kj}y_j\bigg)\bigg\|^2\\
&\geq\beta\sum_{i=1}^mp\sum_{k=1}^p
\bigg\|L_{ki}^*\bigg(T_k\bigg(\sum_{j=1}^mL_{kj}x_j\bigg)
-T_k\bigg(\sum_{j=1}^mL_{kj}y_j\bigg)\bigg)\bigg\|^2\\
&\geq\beta\sum_{i=1}^m \bigg\|\sum_{k=1}^p
L_{ki}^*T_k\bigg(\sum_{j=1}^mL_{kj}x_j\bigg)
-\sum_{k=1}^pL_{ki}^*T_k\bigg(\sum_{j=1}^mL_{kj}y_j\bigg)\bigg\|^2,
\end{aligned}
\end{multline}
which establishes the inequality.
\end{proof}

\begin{remark}
\label{r:25}
Suppose that $T_k\equiv\Id$ in Proposition~\ref{p:coco}.
Then the operators $(B_i)_{1\leq i\leq m}$ of \eqref{e:4oct2007}
are simply those resulting from Proposition~\ref{p:14}. However, since
$\beta_k\equiv 1$, the bound given in \eqref{e:genna07-49} is sharper 
than that given in \eqref{e:genna07-47}.
\end{remark}

\begin{corollary}
\label{c:var}
For every $k\in\{1,\ldots,p\}$, let $\GG_k$ be a real Hilbert space, 
let $\tau_k\in\RPP$, let $\varphi_k\colon\GG_k\to\RR$ be a 
$\tau_k$-Lipschitz-differentiable convex function, and, for every 
$i\in\{1,\ldots,m\}$, let $L_{ki}\in\BL(\HH_i,\GG_k)$. Assume that 
$\min_{1\leq k\leq p}\sum_{i=1}^m\|L_{ki}\|^2>0$ and set 
\begin{equation}
\label{e:13nov2007}
(\forall i\in\{1,\ldots,m\})\quad
B_i\colon\HH_1\times\cdots\times\HH_m\to\HH_i\colon
(x_j)_{1\leq j\leq m}\mapsto\sum_{k=1}^pL_{ki}^*\nabla\varphi_k
\bigg(\sum_{j=1}^mL_{kj}x_j\bigg).
\end{equation}
Then $(B_i)_{1\leq i\leq m}$ satisfies \eqref{e:genna07-21} with 
\begin{equation}
\label{e:genna07-148}
\beta=\frac{1}{p\:\displaystyle{\max_{1\leq k\leq p}}
\:\tau_k\sum_{i=1}^m\|L_{ki}\|^2}.
\end{equation}
\end{corollary}
\begin{proof}
Lemma~\ref{l:BH} asserts that, for every $k\in\{1,\ldots,p\}$, 
$T_k=\nabla\varphi_k$ is $\tau_k^{-1}$-cocoercive. The result 
therefore follows from Proposition~\ref{p:coco}.
\end{proof}

\begin{example}[saddle point problems]
\label{ex:saddle}
For every $k\in\{1,\ldots,p\}$ and $l\in\{1,\ldots,q\}$, let $\GG_k$ 
and ${\mathcal K}_l$ be real Hilbert spaces, let $\tau_k\in\RPP$, 
let $\kappa_l\in\RPP$, let $\varphi_k\colon\GG_k\to\RR$ be a
$\tau_k$-Lipschitz-differentiable convex function,
let $\psi_l\colon{\mathcal K}_l\to\RR$ be a
$\kappa_l$-Lipschitz-differentiable convex function. Furthermore,
for every $i\in\{1,\ldots,m\}$, let ${\mathcal Y}_i$ and 
${\mathcal Z}_i$ be real Hilbert spaces, let 
$F_i\colon{\mathcal Y}_i\oplus {\mathcal Z}_i\to\RXX$ satisfy
\ref{e:saddleK} and \ref{e:saddleK'} in Example~\ref{ex:5}, let
$L_{ki}\in\BL({\mathcal Z}_i,\GG_k)$ and 
$M_{li}\in\BL({\mathcal Y}_i,{\mathcal K}_l)$. It is assumed that 
$\min_{1\leq k\leq p}\sum_{i=1}^m\|L_{ki}\|^2>0$ and that
$\min_{1\leq l\leq q}\sum_{i=1}^m\|M_{li}\|^2>0$. 
Consider the problem 
\begin{equation}
\label{e:prob1CV}
\underset{y_1\in{\mathcal Y}_1,\ldots,y_m\in{\mathcal Y}_m}
{\mathrm{maximize}}\;\;
\underset{z_1\in{\mathcal Z}_1,\ldots,z_m\in{\mathcal Z}_m}
{\mathrm{minimize}}\;\;
\sum_{i=1}^mF_i(y_i,z_i)-\sum_{l=1}^q\psi_l\bigg(\sum_{i=1}^mM_{li}
y_i\bigg)+\sum_{k=1}^p\varphi_k\bigg(\sum_{i=1}^mL_{ki}z_i\bigg).
\end{equation}
Now set 
\begin{equation}
\begin{cases}
\widetilde{B}_{i1}\colon{\mathcal Y}_1\times
\cdots\times{\mathcal Y}_m\to{\mathcal Y}_i
\colon(y_j)_{1\leq j\leq m}\mapsto
\sum_{l=1}^qM_{li}^*\nabla\psi_l
\big(\sum_{j=1}^mM_{lj}y_j\big)\\
\widetilde{B}_{i2}\colon{\mathcal Z}_1\times
\cdots\times{\mathcal Z}_m\to{\mathcal Z}_i
\colon(z_j)_{1\leq j\leq m}\mapsto\sum_{k=1}^pL_{ki}^*
\nabla\varphi_k\big(\sum_{j=1}^mL_{kj}z_j\big),
\end{cases}
\end{equation}
and
\begin{equation}
\label{e:2009-02-10}
\beta_1=\frac{1}{q\max\limits_{1\leq l\leq q}
\kappa_l\sum_{j=1}^m\|M_{lj}\|^2}
\quad\text{and}\quad
\beta_2=\frac{1}{p\max\limits_{1\leq k\leq p}
\tau_k\sum_{j=1}^m\|L_{kj}\|^2}.
\end{equation}
We derive from Corollary~\ref{c:var} that, for every $(y_1,\ldots,y_m)$ 
and $(\overline{y}_1,\ldots,\overline{y}_m)$ in 
${\mathcal Y}_1\times\cdots\times{\mathcal Y}_m$, 
\begin{multline}
\label{e:genna07-216}
\sum_{i=1}^m\scal{\widetilde{B}_{i1}(y_1,\ldots,y_m)-\widetilde{B}_{i1}
(\overline{y}_1,\ldots,\overline{y}_m)}
{y_i-\overline{y}_i}_{{\mathcal Y}_i}\geq\\[-1mm]
\beta_1\sum_{i=1}^m\big\|\widetilde{B}_{i1}(y_1,\ldots,y_m)-
\widetilde{B}_{i1}(\overline{y}_1,\ldots,\overline{y}_m)
\big\|_{{\mathcal Y}_i}^2,
\end{multline}
and that an analogous inequality is satisfied by $\widetilde{B}_{i2}$
with $\beta_2$. On the other hand, using minimax theory \cite{Rock70}, 
we can cast \eqref{e:prob1CV} in the form of \eqref{e:genna07-25o} 
where, for every $i\in\{1,\ldots,m\}$, 
$\HH_i={\mathcal Y}_i\oplus{\mathcal Z}_i$ and
\begin{equation}
\label{e:2009-02-10-}
B_i=(B_{i1},B_{i2})\colon(y_j,z_j)_{1\leq j\leq m}
\mapsto\big(\widetilde{B}_{i1}(y_1,\ldots,y_m),
\widetilde{B}_{i2}(z_1,\ldots,z_m)\big).
\end{equation}
Altogether, it follows from Example~\ref{ex:5} that 
\eqref{e:prob1CV} is a special case of Problem~\ref{prob:64} in which 
$(B_i)_{1\leq i\leq m}$ satisfies \eqref{e:genna07-21} with 
$\beta=\min\{\beta_1,\beta_2\}$.
\end{example}

\section{Coupling evolution inclusions}
\label{sec:6}
Evolution inclusions arise in various fields of applied
mathematics \cite{Hara81,Show97}. In this section, we address the 
problem of solving systems of coupled evolution inclusions with 
periodicity conditions. The notation and definitions introduced 
in Example~\ref{ex:brezis} will be used.

\subsection{Problem formulation and algorithm}

\begin{problem}
\label{prob:69}
Let $({\mathsf H}_i)_{1\leq i\leq m}$ be real Hilbert spaces and let 
$T\in\RPP$. For every $i\in\{1,\ldots,m\}$, set
\begin{equation}
{\mathcal W}_i=\menge{x\in {\EuScript C}([0,T];{\mathsf H}_i)\cap 
W^{1,2}([0,T];{\mathsf H}_i) }{x(T)=x(0)},
\end{equation}
let ${\mathsf f}_i\in\Gamma_0({\mathsf H}_i)$, and let  
${\mathsf B}_{i}\colon{\mathsf H}_1\times\cdots\times
{\mathsf H}_m\to{\mathsf H}_i$. It is assumed that there exists
$\beta\in\RPP$ such that
\begin{multline}
\label{e:genna07-241}
(\forall ({\mathsf x}_1,\ldots,{\mathsf x}_m)\in{\mathsf H}_1
\times\cdots\times{\mathsf H}_m)
(\forall ({\mathsf y}_1,\ldots,{\mathsf y}_m)\in{\mathsf H}_1
\times\cdots\times{\mathsf H}_m)\\
\sum_{i=1}^m\scal{{\mathsf B}_{i}({\mathsf x}_1,\ldots,{\mathsf x}_m)-
{\mathsf B}_{i}({\mathsf y}_1,\ldots,{\mathsf y}_m)}
{{\mathsf x}_i-{\mathsf y}_i}_{{\mathsf H}_i}\geq\beta\sum_{i=1}^m\big
\|{\mathsf B}_{i}({\mathsf x}_1,\ldots,{\mathsf x}_m)-
{\mathsf B}_{i}({\mathsf y}_1,\ldots,{\mathsf y}_m)
\big\|_{{\mathsf H}_i}^2.
\end{multline}
The problem is to 
\begin{multline}
\label{e:genna07-31}
\text{find}\;x_1\in {\mathcal W}_1,\ldots,
x_m\in {\mathcal W}_m\quad\text{such that}\quad \\
(\forall i\in\{1,\ldots,m\})\quad
0\in\displaystyle{x_i'(t)}+
\partial{\mathsf f}_i (x_i(t))+{\mathsf B}_i(x_1 (t),\ldots,x_m (t))\;\;
\text{a.e. on}\;\left]0,T\right[,
\end{multline}
under the assumption that such solutions exist.
\end{problem}

\begin{algorithm}
\label{algo:7}
Fix $\varepsilon\in\left]0,\min\{1,\beta\}\right[$,
$(\gamma_n)_{n\in\NN}$ in $\left[\varepsilon,2\beta-\varepsilon\right]$,
and $(\lambda_n)_{n\in\NN}$ in $\left[0,1-\varepsilon\right]$.
Let, for every $n\in\NN$ and every $i\in\{1,\ldots,m\}$, $y_{i,n}$ be 
the unique solution in ${\mathcal W}_i$ to the inclusion
\begin{multline}
\label{e:11juillet2008-2}
\Frac{x_{i,n}(t)-y_{i,n}(t)}{\gamma_n}
-\big({\mathsf B}_i(x_{1,n}(t),\ldots,x_{m,n}(t))+b_{i,n}(t)\big)\\[2mm]
\in y_{i,n}'(t)+\partial {\mathsf f}_i(y_{i,n}(t))+e_{i,n}(t)
\:\;\text{a.e. on}\;\left]0,T\right[
\end{multline}
and set
\begin{equation}
\label{e:11juillet2008-1}
x_{i,n+1}=\lambda_{i,n}x_{i,n}+(1-\lambda_{i,n})y_{i,n}
\end{equation}
where, for every $i\in\{1,\ldots,m\}$, the following hold.
\begin{enumerate}
\item
\label{c:75}
$x_{i,0}\in W^{1,2}([0,T];{\mathsf H}_i)$.
\item
\label{c:33}
$(b_{i,n})_{n\in\NN}$ and $(e_{i,n})_{n\in\NN}$ are sequences in 
$L^2([0,T];{\mathsf H}_i)$ such that 
\begin{equation}
\label{e:genna07-17}
\sum_{n\in\NN}\sqrt{\int_0^T\|b_{i,n}(t)\|^2_{{\mathsf H}_i}dt}<\pinf
\quad\text{and}\quad 
\sum_{n\in\NN}\sqrt{\int_0^T\|e_{i,n}(t)\|^2_{{\mathsf H}_i}dt}<\pinf.
\end{equation}
\item
\label{c:44}
$(\lambda_{i,n})_{n\in\NN}$ is a sequence in $\left[0,1\right[$
such that $\sum_{n\in\NN}|\lambda_{i,n}-\lambda_n|<\pinf$.
\end{enumerate}
\end{algorithm}

In \eqref{e:11juillet2008-2}, $b_{i,n}(t)$ models the error tolerated in
the computation of ${\mathsf B}_i(x_{1,n}(t),\ldots,x_{m,n}(t))$, while 
$e_{i,n}(t)$ models the error tolerated in solving the inclusion 
with respect to $\partial{\mathsf f}_i(y_{i,n}(t))$. 

\subsection{Convergence}

\begin{theorem}
\label{t:3}
Let $((x_{i,n})_{n\in\NN})_{1\leq i\leq m}$ be sequences generated by
Algorithm~\ref{algo:7}. Then, for every $i\in\{1,\ldots,m\}$, 
$(x_{i,n})_{n\in\NN}$ converges weakly in $W^{1,2}([0,T];{\mathsf H}_i)$ 
to a point $x_i\in {\mathcal W}_i$, and $(x_i)_{1\leq i\leq m}$ is a 
solution to Problem~\ref{prob:69}. 
\end{theorem}
\begin{proof}
For every $i\in\{1,\ldots,m\}$, set $\HH_i=L^2([0,T];{\mathsf H}_i)$ and
\begin{equation}
\label{e:18juillet2008-1}
\begin{array}{lcll}
A_i\colon\!\!\!\!\!&\HH_i&\to&2^{\HH_i}\\[2mm]
&x&\mapsto&
\begin{cases}
\Menge{u\in\HH_i}{u(t)\in x'(t)+\partial {\mathsf f}_i(x(t))\:\;
\text{a.e. in}\;\left]0,T\right[\,},
&\text{if}\;\; x\in {\mathcal W}_i;\\
\emp,&\text{otherwise.}
\end{cases}
\end{array}
\end{equation}
Let us first show that the operators $(A_i)_{1\leq i\leq m}$ are maximal
monotone. For this purpose, let $i\in\{1,\ldots,m\}$, $(x,u)\in\gr A_i$, 
and $(y,v)\in\gr A_i$. It follows from \eqref{e:18juillet2008-1} that, 
almost everywhere on $\left]0,T\right[$,
$u(t)-x'(t)\in\partial {\mathsf f}_i(x(t))$ and
$v(t)-y'(t)\in\partial {\mathsf f}_i(y(t))$. Therefore,
by monotonicity of $\partial {\mathsf f}_i$, we have
\begin{equation}
\label{e:18juillet2008-2}
\int_0^T\scal{\big(u(t)-x'(t)\big)
-\big(v(t)-y'(t)\big)}{x(t)-y(t)}_{{\mathsf H}_i}dt\geq 0.
\end{equation}
Hence,
\begin{align}
\scal{u-v}{x-y}&=\int_0^T\scal{u(t)-v(t)}{x(t)-y(t)}_{{\mathsf H}_i}dt
\nonumber\\
&=\int_0^T\scal{\big(u(t)-x'(t)\big)-\big(v(t)-y'(t)\big)}
{x(t)-y(t)}_{{\mathsf H}_i}dt\nonumber\\
&\quad\;+\int_0^T\scal{x'(t)-y'(t)}
{x(t)-y(t)}_{{\mathsf H}_i}dt\nonumber\\
&\geq\frac12\int_0^T\frac{d\,\|x(t)-y(t)\|
_{{\mathsf H}_i}^2}{dt}dt\nonumber\\
&=\frac12\big(\|x(T)-y(T)\|_{{\mathsf H}_i}^2-\|x(0)-y(0)\|
_{{\mathsf H}_i}^2\big)\nonumber\\
&=0.
\end{align}
Thus, $A_i$ is monotone. To prove maximality, set 
${\mathsf g}_i=(1/2)\|\cdot\|_{{\mathsf H}_i}^2+{\mathsf f}_i$. Then 
${\mathsf g}_i\in\Gamma_0({\mathsf H}_i)$ and 
$\partial{\mathsf g}_i=\Id+\partial {\mathsf f}_i$.
Moreover, since ${\mathsf f}_i\in\Gamma_0({\mathsf H}_i)$, it is 
minorized by a continuous affine functional, say 
${\mathsf f}_i\geq\scal{\cdot}{{\mathsf v}}_{{\mathsf H}_i}+\eta$ for some
${\mathsf v}\in{\mathsf H}_i$ and $\eta\in\RR$.
Now, let ${\mathsf y}\in\dom {\mathsf f}_i=\dom {\mathsf g}_i$ and take 
$({\mathsf x},{\mathsf u})\in\gr\partial {\mathsf g}_i$. Then 
\eqref{e:subdiff} and Cauchy-Schwarz imply the coercivity property
\begin{align}
\frac{\scal{{\mathsf x}-{\mathsf y}}
{{\mathsf u}}_{{\mathsf H}_i}}{\|{\mathsf x}\|_{{\mathsf H}_i}}
&\geq\Frac{{\mathsf g}_i({\mathsf x})-{\mathsf g}_i({\mathsf y})}{
\|{\mathsf x}\|_{{\mathsf H}_i}}\nonumber\\
&=\frac{\|{\mathsf x}\|_{{\mathsf H}_i}}{2}+
\Frac{{\mathsf f}_i({\mathsf x})-
{\mathsf g}_i({\mathsf y})}{\|{\mathsf x}\|_{{\mathsf H}_i}}\nonumber\\
&\geq\frac{\|{\mathsf x}\|_{{\mathsf H}_i}}{2}
-\|{\mathsf v}\|_{{\mathsf H}_i}+
\Frac{\eta-{\mathsf g}_i({\mathsf y})}
{\|{\mathsf x}\|_{{\mathsf H}_i}}\nonumber\\
&\to\pinf\quad\text{as}\quad\|{\mathsf x}\|_{{\mathsf H}_i}\to\pinf. 
\end{align}
Therefore, \cite[Corollaire~3.4]{Brez73} asserts that
for every $w\in\HH_i$ there exists
$z\in{\mathcal W}_i$ such that
\begin{equation}
w(t)\in z'(t)+\partial{\mathsf g}_i(z(t))
=z'(t)+z(t)+\partial{\mathsf f}_i(z(t))\;\:
\text{a.e. on}\;\left]0,T\right[,
\end{equation}
i.e., by \eqref{e:18juillet2008-1}, such that $w-z\in A_iz$. This 
shows that the range of $\Id+A_i$ is $\HH_i$ and hence, by Minty's 
theorem \cite[Theorem~3.5.8]{Aubi90}, that $A_i$ is maximal monotone.

Next, for every $i\in\{1,\ldots,m\}$ and every 
$(x_1,\ldots,x_m)\in\HH_1\times\cdots\times\HH_m$, define 
almost everywhere
\begin{equation}
\label{e:B_i17Juillet}
\begin{array}{lcll}
B_i(x_1,\ldots,x_m)\colon\!\!\!\!\!&[0,T]&\to&{\mathsf H}_i\\[2mm]
&t&\mapsto&{\mathsf B}_i(x_1(t),\ldots,x_m(t)).
\end{array}
\end{equation}
Now let $(x_1,\ldots,x_m)\in\HH_1\times\cdots\times\HH_m$ and
set $(\forall i\in\{1,\ldots,m\})$ 
${\mathsf b}_i={\mathsf B}_i(0,\ldots,0)$. Then it follows from
\eqref{e:genna07-241} and Cauchy-Schwarz that, almost everywhere on
$[0,T]$,
\begin{align}
\beta\sum_{j=1}^m\|{\mathsf B}_j(x_1(t),\ldots,x_m(t))
-{\mathsf b}_j\|_{{\mathsf H}_j}^2
&\leq\sum_{j=1}^m\scal{{\mathsf B}_j(x_1(t),\ldots,x_m(t))
-{\mathsf b}_j}{x_j(t)-0}_{{\mathsf H}_j}\nonumber\\
&\leq\sum_{j=1}^m\|{\mathsf B}_j(x_1(t),\ldots,x_m(t))
-{\mathsf b}_j\|_{{\mathsf H}_j}\,\|x_j(t)\|_{{\mathsf H}_j}\nonumber\\
&\leq\sqrt{\sum_{j=1}^m\|{\mathsf B}_j(x_1(t),\ldots,x_m(t))
-{\mathsf b}_j\|_{{\mathsf H}_j}^2}
\sqrt{\sum_{j=1}^m\|x_j(t)\|_{{\mathsf H}_j}^2}.
\end{align}
Therefore, for every $i\in\{1,\ldots,m\}$,
\begin{align}
\|B_i(x_1,\ldots,x_m)(t)\|_{{\mathsf H}_i}^2
&\leq 2\big(\|{\mathsf b}_i\|_{{\mathsf H}_i}^2+\|B_i(x_1,\ldots,x_m)(t)
-{\mathsf b}_i\|_{{\mathsf H}_i}^2\big)\nonumber\\
&\leq 2\bigg(\|{\mathsf b}_i\|_{{\mathsf H}_i}^2+
\sum_{j=1}^m\|{\mathsf B}_j(x_1(t),\ldots,x_m(t))
-{\mathsf b}_j\|_{{\mathsf H}_j}^2\bigg)\nonumber\\
&\leq2\bigg(\|{\mathsf b}_i\|_{{\mathsf H}_i}^2+
\frac{1}{\beta^2}\sum_{j=1}^m\|x_j(t)\|_{{\mathsf H}_j}^2\bigg)
\:\;\text{a.e. on}\;\left]0,T\right[,
\end{align}
which yields
\begin{equation}
\label{e:santiago2008-09-16}
\int_0^T\|B_i(x_1,\ldots,x_m)(t)\|_{{\mathsf H}_i}^2dt
\leq 2T\|{\mathsf b}_i\|_{{\mathsf H}_i}^2
+\frac{2}{\beta^2}\sum_{j=1}^m\|x_j\|^2,
\end{equation}
so that we can now claim that
$B_i\colon\HH_1\times\cdots\times\HH_m\to L^2([0,T];{\mathsf H}_i)=\HH_i$.
In addition, upon integrating, we derive from \eqref{e:genna07-241} and
\eqref{e:B_i17Juillet} that, for every 
$(y_1,\ldots,y_m)\in\HH_1\times\cdots\times\HH_m$,
\begin{equation}
\label{e:genna07-281}
\sum_{i=1}^m\scal{B_{i}(x_1,\ldots,x_m)-
B_{i}(y_1,\ldots,y_m)}
{x_i-y_i}\geq\beta\sum_{i=1}^m\big
\|B_{i}(x_1,\ldots,x_m)-
B_{i}(y_1,\ldots,y_m)\big\|^2.
\end{equation}
We have thus established \eqref{e:genna07-21}.

Let us now make the connection between Algorithm~\ref{algo:7} and
Algorithm~\ref{algo:1}.  
For every $n\in\NN$ and every $i\in\{1,\ldots,m\}$, 
it follows from \eqref{e:11juillet2008-2}, \eqref{e:18juillet2008-1},
\eqref{e:B_i17Juillet}, and the maximal monotonicity of $A_i$ that 
$y_{i,n}$ is uniquely defined and can be expressed as
\begin{equation}
\label{e:22juillet2008-1}
y_{i,n}=J_{\gamma_n A_i}\Big(x_{i,n}-\gamma_{n}\big
(B_i(x_{1,n},\ldots,x_{m,n})+b_{i,n}\big)\Big)+a_{i,n}\,,
\end{equation}
where
\begin{equation}
\label{e:22juillet2008-2}
a_{i,n}=
J_{\gamma_n A_i}\Big(-\gamma_ne_{i,n}+x_{i,n}-\gamma_{n}\big
(B_i(x_{1,n},\ldots,x_{m,n})+b_{i,n}\big)\Big)-
J_{\gamma_n A_i}\Big(x_{i,n}-\gamma_{n}\big
(B_i(x_{1,n},\ldots,x_{m,n})+b_{i,n}\big)\Big),
\end{equation}
and we therefore derive from
\eqref{e:11juillet2008-2} and \eqref{e:11juillet2008-1} that
\begin{equation}
\label{e:21juillet2008-3}
x_{i,n+1}=\lambda_{i,n}x_{i,n}+(1-\lambda_{i,n})\Big(
J_{\gamma_n A_i}\Big(x_{i,n}-\gamma_{n}\big
(B_i(x_{1,n},\ldots,x_{m,n})+b_{i,n}\big)\Big)+a_{i,n}\Big).
\end{equation}
We observe that \eqref{e:21juillet2008-3} derives from
\eqref{e:genna07-2}, where $A_{i,n}\equiv A_i$ and $B_{i,n}\equiv B_i$.
On the other hand, for every $i\in\{1,\ldots,m\}$,
by nonexpansivity of the operators $(J_{\gamma_n A_i})_{n\in\NN}$,
we deduce from \eqref{e:22juillet2008-2} and \eqref{e:genna07-17} that
\begin{equation}
\sum_{n\in\NN}\|a_{i,n}\|\leq\sum_{n\in\NN}\gamma_n\|e_{i,n}\|
\leq 2\beta\sum_{n\in\NN}\|e_{i,n}\|<\pinf.
\end{equation}
As a result, all the hypotheses of Algorithm~\ref{algo:1} 
are satisfied and hence Theorem~\ref{t:1} asserts that,
for every $i\in\{1,\ldots,m\}$, $(x_{i,n})_{n\in\NN}$ converges weakly 
in $\HH_i=L^2([0,T];{\mathsf H}_i)$ to a point $x_i$, and 
$(x_i)_{1\leq i\leq m}$ satisfies
\begin{equation}
\label{e:genna07-981}
(\forall i\in\{1,\ldots,m\})\quad 0\in A_ix_i+B_{i}(x_1,\ldots,x_m).
\end{equation}
Accordingly, 
\begin{equation}
\label{e:2009-01-22}
\sigma=\max_{1\leq i\leq m}\:\sup_{n\in\NN}\|x_{i,n}\|<\pinf
\end{equation}
and $(\forall i\in\{1,\ldots,m\})$ 
$x_i\in\dom A_i\subset{\mathcal W}_i$. Moreover since, in view of 
\eqref{e:18juillet2008-1} and \eqref{e:B_i17Juillet},
\eqref{e:genna07-981} reduces to \eqref{e:genna07-31}, 
$(x_i)_{1\leq i\leq m}$ is a solution to Problem~\ref{prob:69}. 

To complete the proof, let $i\in\{1,\ldots,m\}$.
To show that $(x_{i,n})_{n\in\NN}$ converges weakly 
to $x_i$ in $W^{1,2}([0,T];{\mathsf H}_i)$, it remains to show
that $(x'_{i,n})_{n\in\NN}$ converges weakly 
to $x_i'$ in $L^2([0,T];{\mathsf H}_i)$. 
We first observe that $(x_{i,n})_{n\in\NN}$ lies in 
$W^{1,2}([0,T];{\mathsf H}_i)$. Indeed, it follows from 
\eqref{e:18juillet2008-1} that
\begin{equation}
(\forall n\in\NN)(\forall z\in\HH_i)\quad 
J_{\gamma_nA_i}z\in\dom(\gamma_n A_i)\subset 
{\mathcal W}_i\subset W^{1,2}([0,T];{\mathsf H}_i). 
\end{equation}
As a result, we deduce from \eqref{e:22juillet2008-2} 
that $(a_{i,n})_{n\in\NN}$ lies in $ W^{1,2}([0,T];{\mathsf H}_i)$.
On the other hand, by construction, $(y_{i,n})_{n\in\NN}$ lies in 
${\mathcal W}_i\subset W^{1,2}([0,T];{\mathsf H}_i)$.
In view of \eqref{e:11juillet2008-1} and~\ref{c:75} in 
Algorithm~\ref{algo:7}, $(x_{i,n})_{n\in\NN}$ 
is therefore in $W^{1,2}([0,T];{\mathsf H}_i)$. Next, let us show that 
$(x_{i,n}')_{n\in\NN}$ is bounded in $L^2([0,T];{\mathsf H}_i)$. 
To this end, let $n\in\NN$ and set
\begin{equation}
\label{e:22juillet2008-3}
w_{i,n}(t)=\Frac{x_{i,n}(t)-y_{i,n}(t)}{\gamma_n}
-{\mathsf B}_i(x_{1,n}(t),\ldots,x_{m,n}(t))
-b_{i,n}(t)-y_{i,n}'(t)-e_{i,n}(t)\:\;\text{a.e. on}\;\left]0,T\right[.
\end{equation}
Then we derive from \eqref{e:11juillet2008-2} that
\begin{equation}
w_{i,n}(t)\in\partial {\mathsf f}_i(y_{i,n}(t))
\:\;\text{a.e. on}\;\left]0,T\right[.
\end{equation}
Hence, since $w_{i,n}\in\HH_i$, it follows from
\cite[Lemme~3.3]{Brez73} that 
\begin{equation}
\label{e:trou}
\frac{d({\mathsf f}_i\circ y_{i,n})(t)}{dt}
=\scal{w_{i,n}(t)}{y_{i,n}'(t)}_{{\mathsf H}_i}
\:\;\text{a.e. on}\;\left]0,T\right[.
\end{equation}
On the other hand, since $y_{i,n}\in{\mathcal W}_i$, we have 
$y_{i,n}(T)=y_{i,n}(0)$. Therefore
\begin{equation}
\label{e:trou2}
\int_0^T\scal{w_{i,n}(t)}{y_{i,n}'(t)}_{{\mathsf H}_i}dt
=\int_0^T \frac{d({\mathsf f}_i\circ y_{i,n})(t)}{dt}dt
={\mathsf f}_i(y_{i,n}(T))-{\mathsf f}_i(y_{i,n}(0))=0
\end{equation}
and, furthermore,
\begin{equation}
\label{e:2009-02-14}
\int_0^T\scal{y_{i,n}(t)}{y_{i,n}'(t)}_{{\mathsf H}_i}dt=
\frac12\int_0^T\frac{d\|y_{i,n}(t)\|_{{\mathsf H}_i}^2}{dt}dt
=\frac{\|y_{i,n}(T)\|_{{\mathsf H}_i}^2
-\|y_{i,n}(0)\|_{{\mathsf H}_i}^2}{2} =0.
\end{equation}
We deduce from \eqref{e:trou2}, \eqref{e:22juillet2008-3}, 
and \eqref{e:2009-02-14} that
\begin{align}
\label{e:22juillet2008-4}
0&
=\int_0^T\Scal{\Frac{x_{i,n}(t)}{\gamma_n}}
{y_{i,n}'(t)}_{{\mathsf H}_i}dt
-\int_0^T\scal{{\mathsf B}_i(x_{1,n}(t),\ldots,x_{m,n}(t))}
{y_{i,n}'(t)}_{{\mathsf H}_i}dt
\nonumber\\[2mm]
&\quad\;
-\int_0^T\scal{b_{i,n}(t)}{y_{i,n}'(t)}_{{\mathsf H}_i}dt
-\int_0^T\|y_{i,n}'(t)\|_{{\mathsf H}_i}^2dt
-\int_0^T\scal{e_{i,n}(t)}{y_{i,n}'(t)}_{{\mathsf H}_i}dt.
\end{align}
Thus, using Cauchy-Schwarz, the inequality $\gamma_n\geq\varepsilon$,
and \eqref{e:B_i17Juillet}, we obtain
\begin{equation}
\label{e:22juillet2008-5}
\|y_{i,n}'\|^2
\leq\Big(\Frac{1}{\varepsilon}\|x_{i,n}\|
+\|B_i(x_{1,n},\ldots,x_{m,n})\|
+\|b_{i,n}\|+\|e_{i,n}\|\Big)\,\|y_{i,n}'\|.
\end{equation}
In turn, it follows from \eqref{e:11juillet2008-1} that
\begin{equation}
\label{e:cmm2008}
\|x_{i,n+1}'\|\leq\lambda_{i,n}\|x_{i,n}'\|+(1-\lambda_{i,n})
\Big(\Frac{1}{\varepsilon}\|x_{i,n}\|+\|B_i(x_{1,n},\ldots,x_{m,n})\|
+\|b_{i,n}\|+\|e_{i,n}\|\Big).
\end{equation}
On the other hand, arguing as in \eqref{e:santiago2008-09-16}, we 
derive from \eqref{e:2009-01-22} that
\begin{equation}
\|B_i(x_{1,n},\ldots,x_{m,n})\|
\leq\sqrt{2T\|{\mathsf b}_i\|_{{\mathsf H}_i}^2
+\frac{2m\sigma^2}{\beta^2}}
\leq
\sqrt{2T}\|{\mathsf b}_i\|_{{\mathsf H}_i}
+\sqrt{2m}\,\frac{\sigma}{\beta}.
\end{equation}
Hence, using \ref{c:33} in Algorithm~\ref{algo:7}, we derive 
by induction from \eqref{e:cmm2008} that
\begin{equation}
\|x_{i,n}'\|\leq\max
\bigg\{\|x_{i,0}'\|,\Frac{\sigma}{\varepsilon}+
\sqrt{2T}\|{\mathsf b}_i\|_{{\mathsf H}_i}+
\sqrt{2m}\,\frac{\sigma}{\beta}
+\sup_{k\in\NN}\big(\|b_{i,k}\|+\|e_{i,k}\|\big)\bigg\}.
\end{equation}
This shows the boundedness of $(x_{i,n}')_{n\in\NN}$ in
$L^2([0,T];{\mathsf H}_i)$. Now let 
$z$ be the weak limit in $L^2([0,T];{\mathsf H}_i)$ of an arbitrary 
weakly convergent subsequence of $(x_{i,n}')_{n\in\NN}$. Since 
$(x_{i,n})_{n\in\NN}$ converges weakly in $L^2([0,T];{\mathsf H}_i)$ to 
$x_i$, it therefore follows from \cite[Proposition~23.19]{Zeid90} that
$z=x_i'$. In turn, this shows that $(x_{i,n}')_{n\in\NN}$
converges weakly in $L^2([0,T];{\mathsf H}_i)$ to $x_i'$.
\end{proof}

\section{The variational case}
\label{sec:5}
In this section, we study a special case of Problem~\ref{prob:64} which
yields a variational formulation that extends \eqref{e:prob5-}.
This framework can be regarded as a particular
instance of Example~\ref{ex:saddle}.

\subsection{Problem formulation and algorithm}
Recall that, for every $f\in\Gamma_0(\HH)$
and every $x\in\HH$, the function $y\mapsto f(y)+\|x-y\|^2/2$ admits a
unique minimizer, which is denoted by $\prox_f x$. The proximity operator
thus defined can be expressed as $\prox_f=J_{\partial f}$ \cite{More65}.

\begin{problem}
\label{prob:62}
Let $(\HH_i)_{1\leq i\leq m}$ and $(\GG_k)_{1\leq k\leq p}$ be real 
Hilbert spaces. For every $i\in\{1,\ldots,m\}$, let 
$f_i\in\Gamma_0(\HH_i)$ and, for every $k\in\{1,\ldots,p\}$, let 
$\tau_k\in\RPP$, let $\varphi_k\colon\GG_k\to\RR$ be a 
$\tau_k$-Lipschitz-differentiable convex function, and let 
$L_{ki}\in\BL(\HH_i,\GG_k)$. It is assumed that 
$\min_{1\leq k\leq p}\sum_{i=1}^m\|L_{ki}\|^2>0$. The problem is to 
\begin{equation}
\label{e:genna07-4}
\underset{x_1\in\HH_1,\ldots,\,x_m\in\HH_m}{\mathrm{minimize}}\;\;
\sum_{i=1}^mf_i(x_i)+
\sum_{k=1}^p\varphi_k\bigg(\sum_{i=1}^mL_{ki}x_i\bigg),
\end{equation}
under the assumption that solutions exist.
\end{problem}

\begin{algorithm}
\label{algo:2}
Set 
\begin{equation}
\label{e:14nov2007}
\beta=\frac{1}{p\:\displaystyle{\max_{1\leq k\leq p}}
\:\tau_k\sum_{i=1}^m\|L_{ki}\|^2}.
\end{equation}
Fix $\varepsilon\in\left]0,\min\{1,\beta\}\right[$,
$(\gamma_n)_{n\in\NN}$ in $\left[\varepsilon,2\beta-\varepsilon\right]$,
$(\lambda_n)_{n\in\NN}$ in $\left[0,1-\varepsilon\right]$, and
$(x_{i,0})_{1\leq i\leq m}$ $\in\HH_1\times\cdots\times\HH_m$.
Set, for every $n\in\NN$,
\begin{equation}
\label{e:genna07-42}
\begin{cases}
x_{1,n+1}=\lambda_{1,n}x_{1,n}\:+\\
\hskip 17mm
(1-\lambda_{1,n})\bigg(\prox_{\gamma_{n}f_{1,n}}\bigg(
x_{1,n}-\gamma_{n}\bigg(\displaystyle{\sum_{k=1}^p}
L_{k1}^*\nabla\varphi_k
\bigg(\sum_{j=1}^mL_{kj}x_{j,n}\bigg)+b_{1,n}\bigg)\bigg)
+a_{1,n}\bigg),\\
~\quad\qquad\:\vdots\\
x_{m,n+1}=\lambda_{m,n}x_{m,n}\:+\\
\hskip 17mm
(1-\lambda_{m,n})\bigg(\prox_{\gamma_{n}f_{m,n}}\bigg(
x_{m,n}-\gamma_{n}\bigg(\displaystyle{\sum_{k=1}^p}
L_{km}^*\nabla\varphi_k
\bigg(\sum_{j=1}^mL_{kj}x_{j,n}\bigg)+b_{m,n}\bigg)\bigg)
+a_{m,n}\bigg),
\end{cases}
\end{equation}
where, for every $i\in\{1,\ldots,m\}$, the following hold.
\begin{enumerate}
\item
\label{c:p1}
$(f_{i,n})_{n\in\NN}$ are functions in $\Gamma_0(\HH_i)$ such that
\begin{equation}
\label{e:p1}
(\forall\rho\in\RPP)\;\;\sum_{n\in\NN}\:\sup_{\|y\|\leq\rho}
\|\prox_{\gamma_n f_{i,n}}y-\prox_{\gamma_nf_i}y\|<\pinf.
\end{equation}
\item
\label{c:p3}
$(a_{i,n})_{n\in\NN}$ and $(b_{i,n})_{n\in\NN}$ are sequences in 
$\HH_i$ such that $\sum_{n\in\NN}\|a_{i,n}\|<\pinf$ and 
$\sum_{n\in\NN}\|b_{i,n}\|<\pinf$.
\item
\label{c:p4}
$(\lambda_{i,n})_{n\in\NN}$ is a sequence in $\left[0,1\right[$
such that $\sum_{n\in\NN}|\lambda_{i,n}-\lambda_n|<\pinf$.
\end{enumerate}
\end{algorithm}

\begin{theorem}
\label{t:2}
Let $((x_{i,n})_{n\in\NN})_{1\leq i\leq m}$ be sequences generated by 
Algorithm~\ref{algo:2}. Then, for every $i\in\{1,\ldots,m\}$, 
$(x_{i,n})_{n\in\NN}$ converges weakly to a point $x_i\in\HH_i$, and 
$(x_i)_{1\leq i\leq m}$ is a solution to Problem~\ref{prob:62}. 
\end{theorem}
\begin{proof}
Problem~\ref{prob:62} is a special case of Problem~\ref{prob:64} where,
for every $i\in\{1,\ldots,m\}$, 
\begin{equation}
\label{e:genna07-5}
A_i=\partial f_i\quad\text{and}\quad
B_i\colon(x_j)_{1\leq j\leq m}\mapsto\sum_{k=1}^p
L_{ki}^*\nabla\varphi_k\bigg(\sum_{j=1}^mL_{kj}x_j\bigg).
\end{equation}
Indeed, define $\HHH$ as in the proof of Theorem~\ref{t:1} and set
\begin{equation}
{\boldsymbol f}\colon\HHH\to\RX\colon(x_i)_{1\leq i\leq m}
\mapsto\sum_{i=1}^mf_i(x_i)
\end{equation}
and
\begin{equation}
{\boldsymbol g}\colon\HHH\to\RR\colon(x_i)_{1\leq i\leq m}\mapsto
\sum_{k=1}^p\varphi_k\bigg(\sum_{i=1}^mL_{ki}x_i\bigg).
\end{equation}
Then ${\boldsymbol f}$ and ${\boldsymbol g}$ are in $\Gamma_0(\HHH)$ and
it follows from Fermat's rule and elementary subdifferential
calculus that, for every $(x_1,\ldots,x_m)\in\HHH$,
\begin{eqnarray}
(x_1,\ldots,x_m)\;\;\text{solves \eqref{e:genna07-4}}\;
&\Leftrightarrow&(0,\ldots,0)\in\partial({\boldsymbol f}+{\boldsymbol g})
(x_1,\ldots,x_m)\nonumber\\
&\Leftrightarrow&(0,\ldots,0)\in\partial{\boldsymbol f}(x_1,\ldots,x_m)
+\nabla{\boldsymbol g}(x_1,\ldots,x_m)\nonumber\\
&\Leftrightarrow&(\forall i\in\{1,\ldots,m\})\;\;0\in\partial f_i(x_i)
+\sum_{k=1}^p
L_{ki}^*\nabla\varphi_k\bigg(\sum_{j=1}^mL_{kj}x_j\bigg)\nonumber\\
&\Leftrightarrow&(\forall i\in\{1,\ldots,m\})\;\;0\in A_ix_i
+B_i(x_1,\ldots,x_m).
\end{eqnarray}
In addition, Lemma~\ref{l:BH} asserts that, for every 
$k\in\{1,\ldots,p\}$, $\nabla\varphi_k$ is $\tau_k^{-1}$-cocoercive. 
In turn, we derive from Corollary~\ref{c:var} that the family 
$(B_i)_{1\leq i\leq m}$ in \eqref{e:genna07-5} satisfies 
\eqref{e:genna07-21} with $\beta$ as in \eqref{e:14nov2007}. Setting
\begin{equation}
(\forall i\in\{1,\ldots,m\})(\forall n\in\NN)\quad
A_{i,n}=\partial f_{i,n}\quad\text{and}\quad
B_{i,n}=B_i,
\end{equation}
we deduce from \eqref{e:p1} that Algorithm~\ref{algo:2} is a particular 
case of Algorithm~\ref{algo:1}. Altogether, Theorem~\ref{t:2} follows 
from Theorem~\ref{t:1}.
\end{proof}

\subsection{Applications}
Let us consider some applications of Theorem~\ref{t:2}, starting with
a game-theoretic interpretation of Problem~\ref{prob:62}.

\begin{example}[coordinated games]
Consider a game with $m$ players indexed by $i\in\{1,\ldots,m\}$. 
The strategy
$x_i$ of the $i$th player lies in the real Hilbert space $\HH_i$ and his
individual utility is modeled by a proper upper semicontinuous concave 
function $h_i\colon\HH_i\to\left[\minf,\pinf\right[$. In the absence of
coordination, the goal of each player is to maximize his own payoff,
which can be described by the variational problem 
\begin{equation}
\label{e:chefd'orchestre1}
\underset{x_1\in\HH_1,\ldots,\,x_m\in\HH_m}{\mathrm{maximize}}\;\;
\sum_{i=1}^mh_i(x_i).
\end{equation}
A coordinator having a global vision of the common interest of the group
of players (say, a benevolent dictator \cite{Moul03}) imposes that,
instead of solving the individualistic problem 
\eqref{e:chefd'orchestre1}, the players solve the joint equilibration 
problem 
\begin{equation}
\label{e:chefd'orchestre2}
\underset{x_1\in\HH_1,\ldots,\,x_m\in\HH_m}{\mathrm{maximize}}\;\;
\sum_{i=1}^mh_i(x_i)+{\boldsymbol g}(x_1,\ldots,x_m),
\end{equation}
where ${\boldsymbol g}\colon\bigoplus_{i=1}^m\HH_i\to\RR$ is a 
Lipschitz-differentiable concave utility function that models the 
collective welfare of the group. A finer model consists in considering 
$p$ subgroups of players and writing 
${\boldsymbol g}=\sum_{k=1}^p{\boldsymbol g}_k$, 
where the payoff ${\boldsymbol g}_k$ of subgroup $k\in\{1,\ldots,p\}$ can
be expressed as
\begin{equation}
\label{e:chefd'orchestre3}
{\boldsymbol g}_k\colon(x_1,\ldots,x_m)\mapsto
\psi_k\bigg(\sum_{i=1}^m L_{ki}x_i\bigg),
\end{equation}
where $\psi_k$ is a Lipschitz-differentiable concave function on a 
real Hilbert space $\GG_k$ and where, for every $i\in\{1,\ldots,m\}$, 
$L_{ki}\in\BL(\HH_i,\GG_k)$. In this model, player $i$ is involved 
in the activity of subgroup $k$ if $L_{ki}\neq 0$. 
Upon setting $f_i=-h_i$ for every $i\in\{1,\ldots,m\}$ and
$\varphi_k=-\psi_k$ for every $k\in\{1,\ldots,p\}$, we recover precisely
Problem~\ref{prob:62}.
Let us notice that a solution $(x_1,\ldots,x_m)$ to 
\eqref{e:chefd'orchestre2}--\eqref{e:chefd'orchestre3} 
can be interpreted as a Nash equilibrium of the potential
game \cite{Mond96} in which the payoff of player $i$ in terms of the 
strategies of the remaining $m-1$ players is given by 
\begin{equation}
\label{e:chefd'orchestre4}
x_i\mapsto h_i(x_i)+\sum_{k=1}^p\psi_k\bigg(\sum_{j=1}^mL_{kj}x_j\bigg).
\end{equation}
In this framework, Theorem~\ref{t:2} provides a numerical construction of 
a Nash equilibrium of the game, and Algorithm~\ref{algo:2} provides
a dynamical model for the interaction between the players. 
At iteration $n$ of Algorithm~\ref{algo:2},
each player $i$ aims at maximizing the utility given in 
\eqref{e:chefd'orchestre4}. This is carried out by the proximal step
\eqref{e:genna07-42}, which is a relaxed version of the exact proximal 
step 
\begin{equation}
\label{e:genna07-42exact}
x_{i,n+1}=\prox_{\gamma_{n}f_{i}}\bigg(
x_{i,n}-\gamma_{n}\sum_{k=1}^pL_{ki}^*\nabla\varphi_k
\bigg(\sum_{j=1}^mL_{kj}x_{j,n}\bigg)\bigg),
\end{equation}
in which the function $f_i$ is replaced by an approximation $f_{i,n}$, 
and some errors $a_{i,n}$ and $b_{i,n}$ are tolerated in the numerical 
implementation of $\prox_{\gamma_{n}f_{i,n}}$ and 
$(\nabla\varphi_k)_{1\leq k\leq p}$,
respectively. The last ingredient of this step concerns risk aversion
and is modeled by the relaxation parameter $\lambda_{i,n}$. When
$\lambda_{i,n}=0$, player $i$ makes a full proximal step; 
otherwise, his step is more heavily anchored to his current position 
$x_{i,n}$ due, for instance,
to uncertainty concerning the next performance of his payoff. Let us note
that, in the absence of coordination ($\varphi_k\equiv 0$) the dynamics of
each player would just evolve independently through pure proximal 
iterations. The coordinator modifies the current strategy $x_{i,n}$
by adding to it a component in the direction of the gradient of 
the collective utility, namely 
$-\gamma_{n}\sum_{k=1}^pL_{ki}^*\nabla\varphi_k
(\sum_{j=1}^mL_{kj}x_{j,n})$.
In this simultaneous dynamical game, the players choose strategies in a
decentralized fashion and without knowledge of the strategies that are 
being chosen by other players. 
\end{example}

\begin{example}[2-agent problem]
In Problem~\ref{prob:62}, set $m=2$ and $p=1$. Then \eqref{e:genna07-4} 
becomes 
\begin{equation}
\label{e:prob4}
\underset{x_1\in\HH_1,\,x_2\in\HH_2}{\mathrm{minimize}}\;\;
f_1(x_1)+f_2(x_2)+\varphi_1(L_{11}x_1+L_{12}x_2).
\end{equation}
Now suppose that $\varphi_1$ is the
Moreau envelope of a function $\psi\in\Gamma_0(\GG_1)$, i.e.,
\begin{equation}
\varphi_1\colon x\mapsto
\;\inf_{y\in\GG_1}\psi(y)+\frac12\|x-y\|^2_{\GG_1}.
\end{equation}
Then $\nabla\varphi_1=\Id-\prox_\psi$ has Lipschitz constant 
$\tau_1=1$ \cite{More65}.
Let us employ the simple form of \eqref{e:genna07-42} in which 
$\lambda_{n}\equiv 0$, $\lambda_{1,n}\equiv 0$, 
$\lambda_{2,n}\equiv 0$, $a_{1,n}\equiv 0$, $a_{2,n}\equiv 0$, 
$f_{1,n}\equiv f_1$, $f_{2,n}\equiv f_2$, 
$b_{1,n}\equiv 0$, and $b_{2,n}\equiv 0$, namely 
\begin{equation}
\label{e:kanshu}
\begin{cases}
x_{1,n+1}=\prox_{\gamma_n f_1}
\big(x_{1,n}+\gamma_nL_{11}^*(\prox_\psi-\Id)
(L_{11}x_{1,n}+L_{12}x_{2,n})\big)\\
x_{2,n+1}=\prox_{\gamma_n f_2}
\big(x_{2,n}+\gamma_nL_{12}^*(\prox_\psi-\Id)
(L_{11}x_{1,n}+L_{12}x_{2,n})\big).
\end{cases}
\end{equation}
Theorem~\ref{t:2} asserts that, if $(\gamma_n)_{n\in\NN}$ lies in
$\left[\varepsilon,2(\|L_{11}\|^2+\|L_{12}\|^2)^{-1}-\varepsilon\right]$ 
for some arbitrarily small $\varepsilon\in\RPP$, then 
$((x_{1,n},x_{2,n}))_{n\in\NN}$ converges weakly to a 
solution $(x_1,x_2)$ to \eqref{e:prob4}. In particular, if 
$\psi$ is the indicator function of a nonempty closed convex subset $C$ 
of $\GG_1$, then \eqref{e:prob4} and \eqref{e:kanshu} become respectively 
\begin{equation}
\label{e:prob7}
\underset{x_1\in\HH_1,\,x_2\in\HH_2}{\mathrm{minimize}}\;\;
f_1(x_1)+f_2(x_2)+\frac12d_C^2(L_{11}x_1+L_{12}x_2)
\end{equation}
and 
\begin{equation}
\label{e:kanshu4}
\begin{cases}
x_{1,n+1}=\prox_{\gamma_n f_1}
\big(x_{1,n}+\gamma_nL_{11}^*(P_C-\Id)
(L_{11}x_{1,n}+L_{12}x_{2,n})\big)\\
x_{2,n+1}=\prox_{\gamma_n f_2}
\big(x_{2,n}+\gamma_nL_{12}^*(P_C-\Id)
(L_{11}x_{1,n}+L_{12}x_{2,n})\big).
\end{cases}
\end{equation}
A further special case of interest is when $C=\{0\}$, meaning that
\eqref{e:prob7} reduces to \eqref{e:prob5-}, i.e.,
\begin{equation}
\label{e:prob5}
\underset{x_1\in\HH_1,\,x_2\in\HH_2}{\mathrm{minimize}}\;\;
f_1(x_1)+f_2(x_2)+\frac12\|L_{11}x_1+L_{12}x_2\|^2_{\GG_1},
\end{equation}
and that \eqref{e:kanshu} assumes the form 
\begin{equation}
\label{e:kanshu2}
\begin{cases}
x_{1,n+1}=\prox_{\gamma_n f_1}
\big(x_{1,n}-\gamma_nL_{11}^*(L_{11}x_{1,n}+L_{12}x_{2,n})\big)\\
x_{2,n+1}=\prox_{\gamma_n f_2}
\big(x_{2,n}-\gamma_nL_{12}^*(L_{11}x_{1,n}+L_{12}x_{2,n})\big).
\end{cases}
\end{equation}
In \cite{Atto08}, \eqref{e:prob5} was approached via an inertial
alternating proximal algorithm. Finally, if we further specialize 
\eqref{e:prob5} by 
choosing $\HH_1=\HH_2=\GG_1$ and $L_{11}=\Id=-L_{12}$, then 
\eqref{e:prob5} reduces 
to \eqref{e:prob62}, which was first considered in \cite{Acke80}.
In this case, upon setting $\gamma_n\equiv 1/2$ in \eqref{e:kanshu2} we 
obtain the parallel proximal algorithm 
\begin{equation}
\label{e:kanshu3}
\begin{cases}
x_{1,n+1}=\prox_{f_1/2}\big((x_{1,n}+x_{2,n})/2\big)\\
x_{2,n+1}=\prox_{f_2/2}\big((x_{1,n}+x_{2,n})/2\big).
\end{cases}
\end{equation}
In view of the above analysis, the sequence 
$((x_{1,n},x_{2,n}))_{n\in\NN}$ thus generated converges 
weakly to a solution to \eqref{e:prob62}. In \cite{Acke80}, the same
conclusion was reached for the sequential algorithm (see 
also \cite{Atto07} for an alternative algorithm with costs-to-move)
\begin{equation}
\label{e:kanshu5}
\begin{cases}
x_{1,n+1}=\prox_{f_1}x_{2,n}\\
x_{2,n+1}=\prox_{f_2}x_{1,n+1}.
\end{cases}
\end{equation}
\end{example}

\begin{example}[traffic theory]
\label{ex:transport}
Consider a network with $M$ links indexed by $j\in\{1,\ldots,M\}$ and 
$N$ paths indexed by $l\in\{1,\ldots,N\}$, linking a subset of $Q$
origin-destination node pairs indexed by $k\in\{1,\ldots,Q\}$. 
There are $m$ types of users 
indexed by $i\in\{1,\ldots,m\}$ transiting on the network. For every 
$i\in\{1,\ldots,m\}$ and $l\in\{1,\ldots,N\}$, let
$\xi_{il}\in\RR$ be the flux of user $i$ on path $l$ and 
let $x_i=(\xi_{il})_{1\leq l\leq N}$ be the flow associated with
user $i$. A standard problem in traffic theory is to find 
a Wardrop equilibrium \cite{Ward52} of the network, i.e.,
flows $(x_i)_{1\leq i\leq m}$ such that the costs in all paths actually
used are equal and less than those a single user would face on 
any unused path. Such an equilibrium can be obtained by solving the 
variational problem \cite{Beck56,Patr94,Shef85}
\begin{equation}
\label{e:Beck}
\underset{x_1\in C_1,\ldots,\,x_m\in C_m}{\mathrm{minimize}}\;\;
\sum_{j=1}^M\int_0^{h_j(x_1,\ldots,x_m)}\phi_j(h)dh,
\end{equation}
where $\phi_j\colon\RR\to\RP$ is a strictly increasing $\tau$-Lipschitz 
continuous function modeling the cost of transiting on link $j$ and 
$h_j(x_1,\ldots,x_m)$ is the total flow through link $j$, which can
be expressed as $h_j(x_1,\ldots,x_m)=\sum_{i=1}^m(Lx_i)^{\top}e_j$,
where $e_j$ is the $j$th canonical basis vector of $\RR^M$ and $L$ 
is an $M\times N$ binary matrix with $jl$th entry equal to $1$ 
or $0$, according as link $j$ belongs to path $l$ or not. Furthermore, 
each closed and convex constraint set $C_i$ in 
\eqref{e:Beck} is defined 
as $C_i=\menge{(\eta_l)_{1\leq l\leq N}\in\RP^N}
{(\forall k\in \{1,\ldots, Q\})\;\sum_{l\in N_k}\eta_l=\delta_{ik}}$,
where $\emp\neq N_k\subset\{1,\ldots,N\}$ is the set of paths 
linking the pair $k$ and $\delta_{ik}\in\RP$ is the flow of user 
$i$ that must transit from the origin to the destination 
of pair $k$ (for more details on network flows, see 
\cite{Rock84,Shef85}). Upon setting
\begin{equation}
\label{e:var1}
\varphi_1\colon\RR^M\to\RR\colon(\nu_j)_{1\leq j\leq M}\mapsto
\sum_{j=1}^M\int_0^{\nu_j}\phi_j(h)dh,
\end{equation} 
problem \eqref{e:Beck} can be written as
\begin{equation}
\label{e:probT1}
\underset{x_1\in\RR^N,\ldots,\,x_m\in\RR^N}{\mathrm{minimize}}\;\;
\sum_{i=1}^m\iota_{C_i}(x_i)+\varphi_1\bigg(\sum_{i=1}^mLx_i\bigg).
\end{equation}
Since $\varphi_1$ is strictly convex and
$\tau$-Lipschitz-differentiable,
\eqref{e:probT1} is a particular instance of Problem~\ref{prob:62} with 
$p=1$, $\GG_1=\RR^M$ and $(\forall i\in\{1,\ldots,m\})$
$\HH_i=\RR^N$, $f_i=\iota_{C_i}$, and $L_{1i}=L$.
Accordingly, Theorem~\ref{t:2} asserts that \eqref{e:probT1} can be 
solved by Algorithm~\ref{algo:2} which, with the choice of parameters 
$\gamma_n\equiv\gamma\in\left]0,2/\tau\right[$,
$\lambda_{i,n}\equiv 0$, $\lambda_{n}\equiv 0$, 
$a_{i,n}\equiv 0$, and $b_{i,n}\equiv 0$, yields
\begin{equation}
\label{e:dimitri}
(\forall i\in\{1,\dots,m\})\quad 
x_{i,n+1}=P_{C_i}\bigg(x_{i,n}-\gamma L^\top\bigg(
\phi_1\bigg(\sum_{i=1}^mLx_{i,n}\bigg),\ldots,
\phi_m\bigg(\sum_{i=1}^mLx_{i,n}\bigg)\bigg)\bigg).
\end{equation}
In the special case when $m=1$ the algorithm described in 
\eqref{e:dimitri} is proposed in \cite{Bert82}. Let us note that,
as an alternative to $\varphi_1$ in \eqref{e:var1}, we can consider 
the function
\begin{equation}
\varphi_1\colon\RR^M\to\RR\colon(\nu_j)_{1\leq j\leq M}
\mapsto\sum_{j=1}^M\nu_j\phi_j(\nu_j),
\end{equation}
under suitable assumptions on $(\phi_j)_{1\leq j\leq M}$. In this case, 
\eqref{e:probT1} reduces to the problem of finding the social 
optimum in the network \cite{Shef85}, that is 
\begin{equation}
\underset{x_1\in C_1,\ldots,\,x_m\in C_m}{\mathrm{minimize}}\;\;
\sum_{j=1}^Mh_j(x_1,\ldots,x_m)\phi_j\big(h_j(x_1,\ldots,x_m)\big),
\end{equation}
which can also be solved with Algorithm~\ref{algo:2}.
\end{example}

\begin{example}[source separation]
\label{ex:sign1}
Consider the problem of recovering $m$ signals $(x_i)_{1\leq i\leq m}$
lying in respective Hilbert spaces $(\HH_i)_{1\leq i\leq m}$ from $p$
observations $(z_k)_{1\leq k\leq p}$ lying in respective Hilbert spaces
$(\GG_k)_{1\leq k\leq p}$. The data formation model is
\begin{equation}
\label{e:obs}
(\forall k\in\{1,\ldots,p\})\quad z_k=\sum_{i=1}^mL_{ki}x_i+w_k,
\end{equation}
where $L_{ki}\in\BL(\HH_i,\GG_k)$ and where $w_k\in\GG_k$ models
observation noise (see in particular \cite{Bobi07,Joho01}). 
In other words, the objective is to recover the original signals 
$(x_i)_{1\leq i\leq m}$ from the $p$ mixtures $(z_k)_{1\leq k\leq p}$.
This situation arises in particular in audio signal processing, when
$p$ microphones record the superpositions $(z_k)_{1\leq k\leq p}$ 
of $m$ sources $(x_i)_{1\leq i\leq m}$ that have undergone linear
distortions and noise corruption. 
Let us note that the same type of model arises in
multicomponent signal deconvolution problems \cite{Anth06,Kang97}. 
A variational formulation of the problem is
\begin{equation}
\label{e:cardoso}
\underset{x_1\in\HH_1,\ldots,\,x_m\in\HH_m}{\mathrm{minimize}}\;\;
\sum_{i=1}^mf_i(x_i)+\sum_{k=1}^pD_k\bigg(\sum_{i=1}^mL_{ki}x_i,z_k\bigg).
\end{equation}
In this formulation, each function $f_i\in\Gamma_0(\HH_i)$ models some 
prior knowledge about the signal $x_i$. On the other hand, each function 
$D_k\colon\GG_k\times\GG_k\to\RP$ promotes data fitting: it vanishes 
only on the diagonal $\menge{(z,z)}{z\in\GG_k}$ and, for every $z\in\GG_k$,
$D_k(\cdot,z)$ is convex and Lipschitz-differentiable (for instance, 
$D_k$ can be a Bregman distance under suitable 
assumptions \cite{Sico03,Cens97}, and in particular the standard
quadratic fitting term $D_k\colon(y,z)\mapsto\|y-z\|_{\GG_k}^2$). 
It is clear that \eqref{e:cardoso} 
is a special realization of Problem~\ref{prob:62} (with 
$\varphi_k=D_k(\cdot,z_k)$ for every $k\in\{1,\ldots,p\}$) and that it 
can therefore be solved via Algorithm~\ref{algo:2}.
\end{example}

\begin{example}[image decomposition]
\label{ex:sign2}
A standard problem in image processing is to find the decomposition 
$(x_i)_{1\leq i\leq m}$ of an image $x=\sum_{i=1}^mx_i$ in some 
Hilbert space $\HH$, from some observation $z$. 
When $m=2$, a common instance of this problem is the geometry/texture 
decomposition problem \cite{Aujo04,Aujo06}. The variational 
formulations studied in these papers are special instances of the 
problem 
\begin{equation}
\label{e:prob05}
\underset{x_1\in\HH,\ldots,\,x_m\in\HH}{\mathrm{minimize}}\;\;
\sum_{i=1}^mf_i(x_i)+\frac{1}{4}\bigg\|z-\sum_{i=1}^mx_i\bigg\|^2,
\end{equation}
where $(f_i)_{1\leq i\leq m}$ are functions in $\Gamma_0(\HH)$. The first 
term in the objective is a separable function, the purpose of which is to
promote certain known features of each component $x_i$, 
and the second is a least-squares data fitting term. 
As shown in \cite{Smms05}, for $m=2$, \eqref{e:prob05} can 
be solved by alternating proximal methods, which produce weakly 
convergent sequences. In \cite{Aujo05}, a finer 3-component model 
of the form \eqref{e:prob05} was investigated in $\HH=\RR^N$, 
and a coordinate descent 
algorithm was proposed to solve it. This algorithm, however, has modest 
convergence properties, and it was proved only that the cluster points 
of the sequence it generates are solutions of the particular finite 
dimensional problem considered there. 
By contrast, since \eqref{e:prob05} is a special case of
Problem~\ref{prob:62} (with $\HH_i\equiv\HH$, $k=1$,
$\varphi_1=\|z-\cdot\|^2/4$, and $L_{1i}\equiv\Id$), we can derive 
from Theorem~\ref{t:2} an iterative method 
the orbits of which are guaranteed to converge weakly to a solution to 
\eqref{e:prob05}, under the sole assumption that solutions exist. 
For instance, for $m=3$, \eqref{e:14nov2007} yields $\beta=2/3$. 
Taking for simplicity $\gamma_n\equiv 1$, 
$\lambda_{n}\equiv 0$, and, for $i\in\{1,2,3\}$, $\lambda_{i,n}\equiv 0$, 
$a_{i,n}\equiv 0$, $f_{i,n}\equiv f_i$, and $b_{i,n}\equiv 0$, 
\eqref{e:genna07-42} becomes
\begin{equation}
\label{e:seattle08}
\begin{cases}
x_{1,n+1}=\prox_{f_1}\big((z+x_{1,n}-x_{2,n}-x_{3,n})/2\big)\\
x_{2,n+1}=\prox_{f_2}\big((z-x_{1,n}+x_{2,n}-x_{3,n})/2\big)\\
x_{3,n+1}=\prox_{f_3}\big((z-x_{1,n}-x_{2,n}+x_{3,n})/2\big).
\end{cases}
\end{equation}
Let us note that, since Theorem~\ref{t:2} allows for more general 
coupling terms than that used in \eqref{e:prob05}, more sophisticated 
image decomposition problems can be solved in our framework. 
\end{example}

\begin{example}[best approximation]
\label{ex:best approx}
The convex feasibility problem is to find a point in the intersection of
closed convex subsets $(C_i)_{1\leq i\leq m}$ of a real Hilbert 
space $\HH$ \cite{Baus96,Proc93}. In many instances, this intersection 
may turn out to be 
empty and a relaxation of this problem in the presence of a hard 
constraint represented by $C_1$ is to \cite{Sign99}
\begin{equation}
\label{e:bestapp1}
\underset{x_1\in C_1}{\mathrm{minimize}}\;\;
\frac12\sum_{i=2}^{m}\omega_{i}d_{C_{i}}^2(x_1),
\end{equation}
where $(\omega_i)_{2\leq i\leq m}$ are strictly positive weights
such that $\max_{2\leq i\leq m}\omega_i=1$.
Since, for every $i\in\{2,\ldots,m\}$ and every $x_1\in C_1$,
$d_{C_{i}}^2(x_1)=\min_{x_{i}\in C_{i}}\|x_1-x_{i}\|^2$, 
\eqref{e:bestapp1} can be reformulated as 
\begin{equation}
\label{e:bestapp}
\underset{x_1\in C_1,\ldots,\,x_m\in C_m}{\mathrm{minimize}}\;\;
\frac12\sum_{k=1}^{m-1}\omega_{k+1}\|x_1-x_{k+1}\|^2.
\end{equation}
This is a special instance of Problem~\ref{prob:62} with
$p=m-1$ and, for every $i\in\{1,\ldots,m\}$, $f_i=\iota_{C_i}$ and
\begin{equation}
(\forall k\in\{1,\ldots,m-1\})\quad
\varphi_k=\frac{\omega_{k+1}}{2}\,\|\cdot\|^2\quad\text{and}\quad L_{ki}=
\begin{cases}
\Id,&\text{if}\;\;i=1;\\
-\Id,&\text{if}\;\;i=k+1;\\
0,&\text{otherwise.}
\end{cases}
\end{equation}
We can derive from Algorithm~\ref{algo:2} an algorithm which,
by Theorem~\ref{t:2}, generates orbits that are guaranteed to 
converge weakly to a solution to \eqref{e:bestapp}. Indeed, in this 
case, \eqref{e:14nov2007} yields $\beta=1/(2(m-1))$. For example,
upon setting $\gamma_{n}\equiv\gamma\in\left]0,1/(m-1)\right[$, 
$\lambda_{n}\equiv 0$, $\lambda_{i,n}\equiv 0$, 
$a_{i,n}\equiv 0$, $b_{i,n}\equiv 0$, and $f_{i,n}=\iota_{C_i}$ 
for simplicity, Algorithm~\ref{algo:2} becomes 
\begin{equation}
\label{e:algo3}
\begin{cases}
x_{1,n+1}=P_{C_1}\big((1-\gamma\sum_{i=2}^m\omega_i)x_{1,n}+
\gamma\sum_{i=2}^m\omega_ix_{i,n}\big)\\
(\forall i\in\{2,\ldots,m\})\quad 
x_{i,n+1}=P_{C_i}\big(\gamma\omega_ix_{1,n}+
(1-\gamma\omega_i)x_{i,n}\big).
\end{cases}
\end{equation}
In the particular case when $m=2$ and $\gamma=1/2$,
then $\omega_2=1$, \eqref{e:bestapp} is equivalent to finding a 
best approximation pair relative to $(C_1,C_2)$ \cite{Baus04},
and \eqref{e:algo3} reduces to
\begin{equation}
\begin{cases}
x_{1,n+1}=P_{C_1}\big((x_{1,n}+x_{2,n})/2\big)\\
x_{2,n+1}=P_{C_2}\big((x_{1,n}+x_{2,n})/2\big).
\end{cases}
\end{equation}
\end{example}

\section{Variational problems over decomposed domains in Sobolev spaces}
\label{sec:7}
In this section, we consider a particular case of Problem~\ref{prob:62}
involving Sobolev trace operators in coupling terms modeling constraints
or transmission conditions at the interfaces of subdomains. 

\begin{figure}[!ht]
\begin{center}
\setlength{\unitlength}{0.0176mm}
\begingroup\makeatletter\ifx\SetFigFont\undefined%
\gdef\SetFigFont#1#2#3#4#5{%
  \reset@font\fontsize{#1}{#2pt}%
  \fontfamily{#3}\fontseries{#4}\fontshape{#5}%
  \selectfont}%
\fi\endgroup%
{\renewcommand{\dashlinestretch}{30}
\begin{picture}(6900,3127)(0,-10)
\thicklines
\put(793.063,4700.694){\arc{6785.184}{1.3506}{1.8028}}
\put(36.492,1102.909){\arc{3069.366}{5.3061}{6.8330}}
\put(2207.000,3323.250){\arc{4569.684}{1.2883}{1.8533}}
\put(868.090,1233.208){\arc{3957.029}{5.4206}{6.9104}}
\put(2423.509,1231.296){\arc{4305.230}{5.4461}{6.8624}}
\put(3996.133,1258.798){\arc{3738.473}{5.4549}{6.8013}}
\put(4794.000,-3600.250){\arc{10733.752}{4.6591}{4.9034}}
\put(3450,1425){\ellipse{6884}{2834}}
\put(630,880){\makebox(0,0)[lb]{\smash{{\SetFigFont{11}{14.4}%
{\rmdefault}{\mddefault}{\updefault}$\Omega_2$}}}}
\put(2100,1700){\makebox(0,0)[lb]{\smash{{\SetFigFont{11}{14.4}%
{\rmdefault}{\mddefault}{\updefault}$\Omega_3$}}}}
\put(1900,430){\makebox(0,0)[lb]{\smash{{\SetFigFont{11}{14.4}%
{\rmdefault}{\mddefault}{\updefault}$\Omega_4$}}}}
\put(4650,2270){\makebox(0,0)[lb]{\smash{{\SetFigFont{11}{14.4}%
{\rmdefault}{\mddefault}{\updefault}$\Omega_{m-2}$}}}}
\put(6270,1444){\makebox(0,0)[lb]{\smash{{\SetFigFont{11}{14.4}%
{\rmdefault}{\mddefault}{\updefault}$\Omega_m$}}}}
\put(700,1850){\makebox(0,0)[lb]{\smash{{\SetFigFont{11}{14.4}%
{\rmdefault}{\mddefault}{\updefault}$\Omega_{1}$}}}}
\put(40,2150){\makebox(0,0)[lb]{\smash{{\SetFigFont{10}{14.4}%
{\rmdefault}{\mddefault}{\updefault}$\Upsilon_{1,1}$}}}}
\put(400,450){\makebox(0,0)[lb]{\smash{{\SetFigFont{10}{14.4}%
{\rmdefault}{\mddefault}{\updefault}$\Upsilon_{2,2}$}}}}
\put(6400,530){\makebox(0,0)[lb]{\smash{{\SetFigFont{10}{14.4}%
{\rmdefault}{\mddefault}{\updefault}$\Upsilon_{m,m}$}}}}
\put(3505,1309){\makebox(0,0)[lb]{\smash{{\SetFigFont{17}{20.4}%
{\rmdefault}{\mddefault}{\updefault}. . .}}}}
\put(5710,2029){\makebox(0,0)[lb]{\smash{{\SetFigFont{10}{12.0}%
{\rmdefault}{\mddefault}{\updefault}$\Upsilon_{m-2,m}$}}}}
\put(4800,1824){\makebox(0,0)[lb]{\smash{{\SetFigFont{10}{12.0}%
{\rmdefault}{\mddefault}{\updefault}$\Upsilon_{m-1,m-2}$}}}}
\put(4930,994){\makebox(0,0)[lb]{\smash{{\SetFigFont{11}{14.4}%
{\rmdefault}{\mddefault}{\updefault}$\Omega_{m-1}$}}}}
\put(5890,994){\makebox(0,0)[lb]{\smash{{\SetFigFont{10}{12.0}%
{\rmdefault}{\mddefault}{\updefault}$\Upsilon_{m-1,m}$}}}}
\put(1350,1910){\makebox(0,0)[lb]{\smash{{\SetFigFont{10}{12.0}%
{\rmdefault}{\mddefault}{\updefault}$\Upsilon_{1,3}$}}}}
\put(1570,1219){\makebox(0,0)[lb]{\smash{{\SetFigFont{10}{12.0}%
{\rmdefault}{\mddefault}{\updefault}$\Upsilon_{2,3}$}}}}
\put(670,1394){\makebox(0,0)[lb]{\smash{{\SetFigFont{10}{12.0}%
{\rmdefault}{\mddefault}{\updefault}$\Upsilon_{1,2}$}}}}
\put(2110,870){\makebox(0,0)[lb]{\smash{{\SetFigFont{10}{12.0}%
{\rmdefault}{\mddefault}{\updefault}$\Upsilon_{3,4}$}}}}
\put(1100,634){\makebox(0,0)[lb]{\smash{{\SetFigFont{10}{12.0}%
{\rmdefault}{\mddefault}{\updefault}$\Upsilon_{2,4}$}}}}
\put(3475,1829){\makebox(0,0)[lb]{\smash{{\SetFigFont{11}{14.4}%
{\rmdefault}{\mddefault}{\updefault}$\Omega$}}}}
\end{picture}
}
\caption{Decomposition of the domain $\Omega$.}
\label{fig:example}
\end{center}
\end{figure}
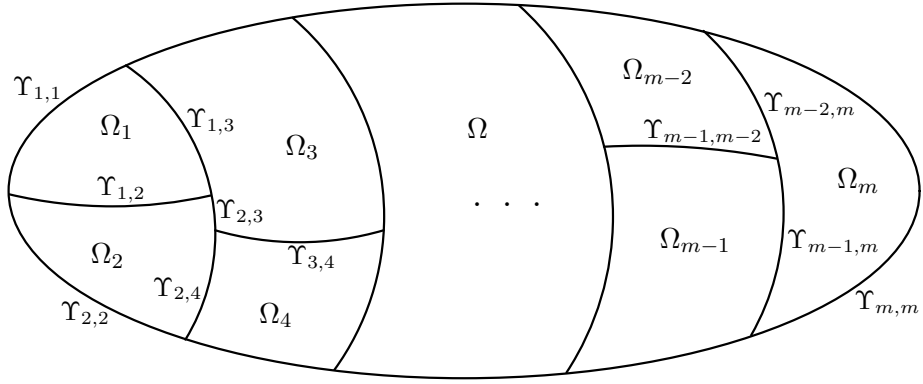

\subsection{Notation and definitions}

We set some notation and recall basic definitions. For details and
complements, see \cite{Adam03,Drab07,Gris85,Neca67,Zeid90}.

We denote by $\RR^N$ the usual $N$-dimensional Euclidean space and by 
$|\cdot|$ its norm, where $N\geq 2$. Let $\Omega$ be a nonempty open 
bounded subset of $\RR^N$ with Lipschitz boundary $\bdry\Omega$.  
The space $H^1(\Omega)=\menge{x\in L^2(\Omega)}{Dx\in(L^2(\Omega))^N}$, 
where $D$ denotes the weak gradient, is a Hilbert space with scalar 
product $\scal{\cdot}{\cdot}_{H^1(\Omega)}\colon(x,y)\mapsto
\int_{\Omega}xy+\int_{\Omega}(Dx)^\top Dy$. We denote by $S$ the surface 
measure on $\bdry\Omega$ \cite[Section~1.1.3]{Neca67}. Now let 
$\Upsilon$ be a nonempty open set 
in $\bdry\Omega$ 
and let $L^2(\Upsilon)$ be the space of square $S$-integrable 
functions on $\Upsilon$. Endowed with the scalar product 
\begin{equation}
\label{e:L2Ups}
\scal{\cdot}{\cdot}_{L^2(\Upsilon)}\colon(v,w)
\mapsto\int_{\Upsilon}vw\,dS,
\end{equation}
$L^2(\Upsilon)$ is a Hilbert space.
The Sobolev trace operator associated with $\Omega$ is the
unique operator
${\mathsf T}\in\BL(H^1(\Omega),L^2(\bdry\Omega))$
such that $(\forall x\in{\EuScript C}^1(\overline{\Omega}))$
${\mathsf T}x=x|_{\bdry\Omega}$. 
Endowed with the scalar product 
\begin{equation}
\label{e:scaledp}
\scal{\cdot}{\cdot}\colon(x,y)\mapsto\int_{\Omega}(D x)^\top D y,
\end{equation}
the space $H^1_{0,\Upsilon}(\Omega)=\menge{x\in H^1(\Omega)}
{{\mathsf T}\,x=0\:\;\text{on}\:\Upsilon}$
is a Hilbert space \cite[Section~25.10]{Zeid90}. Finally,
for $S$-almost every $\omega\in\bdry\Omega$, there exists a unit outward
normal vector $\nu(\omega)$. 

\subsection{Problem formulation and algorithm}

\begin{problem}
\label{prob:edp}
Let $\Omega$ be a nonempty open bounded subset of $\RR^N$ with Lipschitz 
boundary $\bdry\Omega$. Let $(\Omega_i)_{1\leq i\leq m}$ be disjoint 
open subsets of $\Omega$ (see Fig.~\ref{fig:example}) such that the 
boundaries $(\bdry\Omega_i)_{1\leq i\leq m}$ are Lipschitz, 
$\overline{\Omega}=\bigcup_{i=1}^m\overline{\Omega_i}$, and
\begin{equation}
\label{e:15octobre2008}
(\forall i\in\{1,\ldots,m\})\quad\Upsilon_{i,i}=
\inte_{\bdry\Omega}(\bdry\Omega_i\cap\bdry\Omega)\neq\emp,
\end{equation}
where $\inte_{\bdry\Omega}$ denotes the interior relative to
${\bdry\Omega}$. For every $i\in\{1,\ldots,m\}$, set 
\begin{equation}
\label{e:24fevrier2009}
(\forall j\in\{i+1,\ldots,m\})\quad\Upsilon_{i,j}
=\Upsilon_{j,i}=\inte_{\bdry\Omega_i}(\bdry\Omega_i\cap\bdry\Omega_j), 
\end{equation}
let 
\begin{equation}
J(i)=\menge{j\in\{1,\ldots,m\}
\smallsetminus\{i\}}{\Upsilon_{i,j}\neq\emp}, 
\end{equation}
be the set of indices of active interfaces,
let ${\mathsf T}_{i}\colon H^1(\Omega_i)\to L^2(\bdry\Omega_{i})$ be 
the trace operator, let 
\begin{equation}
\label{e:HH_iEDP}
\HH_i=H^1_{0,\Upsilon_{i,i}}(\Omega_i)=\menge{x\in H^1(\Omega_i)}
{{\mathsf T}_i\,x=0\:\;\text{on}\:\Upsilon_{i,i}},
\end{equation}
let $f_i\in\Gamma_0(\HH_i)$, and, for every $j\in J(i)$, 
let $\tau_{ij}\in\RPP$, let
$\varphi_{ij}\colon L^2(\Upsilon_{i,j})\to\RR$ be convex and
$\tau_{ij}$-Lipschitz-differentiable, and set
${\mathsf T}_{ij}\colon\HH_i\to L^2(\Upsilon_{i,j})\colon x\mapsto
({\mathsf T}_ix)|_{\Upsilon_{i,j}}$.
The problem is to
\begin{equation}
\label{prob:16}
\underset{x_1\in \HH_1,\ldots,x_m\in \HH_m}{\mathrm{minimize}}\;\;
\sum_{i=1}^mf_i(x_i)+
\sum_{i=1}^{m}\sum_{j\in J(i)}\varphi_{ij} 
({\mathsf T}_{ij}\,x_i-{\mathsf T}_{ji}\,x_j),
\end{equation}
under the assumption that solutions exist.
\end{problem}

In the above formulation, each function $x_i$ is defined on a
subdomain $\Omega_i$. The potential $f_i$ models intrinsic properties of
$x_i$ while, for every $j\in J(i)$, the potential $\varphi_{ij}$ arising 
in the coupling term models the interaction with the $j$th subdomain as 
a function of the difference of the Sobolev traces of $x_i$ and 
$x_j$ on $\Upsilon_{i,j}$, i.e., of the jump across the interface 
between $\Omega_i$ and $\Omega_j$.
Such variational formulations arise in the modeling of transmission
problems through thin layers, of Neumann's sieve (transmission through 
a finely perforated surface), and of cracks in material
\cite{Atto84,Atto08,Siam06,Atto85}. 
Note that, contrary to these approaches,
our setting can handle $m>2$ domains as well as nonquadratic 
functions $\varphi_{ij}$. We also observe that, if each 
$\varphi_{ij}\colon L^2(\Upsilon_{i,j})\to\RP$ and
vanishes only at $0$, \eqref{prob:16} can be
regarded as a relaxation of some domain decomposition problems, in which
one typically imposes the ``no-jump" conditions 
${\mathsf T}_{ij}\,x_i={\mathsf T}_{ji}\,x_j$ across interfaces
\cite{Bour88,Quar99,Tose05}. More generally, \eqref{prob:16} 
can promote various
properties of the jump. For instance, if $\varphi_{ij}=d_{C_{ij}}^2$,
where $C_{ij}$ is a closed convex subset of $L^2(\Upsilon_{i,j})$, 
the underlying constraint is
${\mathsf T}_{ij}\,x_i-{\mathsf T}_{ji}\,x_j\in C_{ij}$. Unilateral
conditions \cite{Daml85,Khlu04} can be modeled in this fashion.

\begin{algorithm}
\label{algo:edp}
Set
\begin{equation}
\label{e:edp16oct2008}
\beta=\frac{1}{p\:\max\limits_{(k,l)\in\KK}
\tau_{kl}\big(\|{\mathsf T}_{kl}\|^2+\|{\mathsf T}_{lk}\|^2\big)},
\end{equation}
where $p$ is the cardinality of 
$\KK=\menge{(k,l)}{1\leq k\leq m,\;l\in J(k)}$,
and fix $\varepsilon\in\left]0,\min\{1,\beta\}\right[$, 
$(\gamma_n)_{n\in\NN}$ in 
$\left[\varepsilon,2\beta-\varepsilon\right]$, and 
$(\lambda_n)_{n\in\NN}$ in $\left[0,1-\varepsilon\right]$.
For every $n\in\NN$ and every $i\in\{1,\ldots,m\}$, let $y_{i,n}$ be 
the unique solution in $\HH_i$ to the problem
\begin{equation}
\label{e:edp0}
\underset{y\in\HH_i}{\mathrm{minimize}}\;\;
\gamma_nf_i(y)
+\frac{1}{2}\int_{\Omega_i}|D y-D x_{i,n}
+\gamma_n D (z_{i,n}+ b_{i,n})|^2,
\end{equation}
where $z_{i,n}$ is the unique weak solution in $H^1(\Omega_i)$ to the 
Dirichlet-Neumann boundary problem ($\nu_i(\omega)$ is the unit 
outward normal vector at $\omega\in\bdry\Omega_i$) 
\begin{equation}
\label{e:edp1}
\begin{cases}
\Delta z_{i,n}=0&\text{on}\;\Omega_i\\
z_{i,n}=0&\text{on}\;\Upsilon_{i,i}\\
\nu_i^\top D z_{i,n}={\displaystyle\sum_{j\in J(i)}}\widetilde{v_{ij,n}}
&\text{on}\;{\displaystyle\bigcup_{j\in J(i)}}\Upsilon_{i,j},
\end{cases}
\end{equation}
where, for every $j\in J(i)$,
\begin{equation}
\label{e:25fevrier2009}
\widetilde{v_{ij,n}}=
\begin{cases}
v_{ij,n}=
\nabla\varphi_{ij}({\mathsf T}_{ij}\,x_{i,n}-
{\mathsf T}_{ji}\,x_{j,n})-
\nabla\varphi_{ji}({\mathsf T}_{ji}\,x_{j,n}-
{\mathsf T}_{ij}\,x_{i,n})&\text{on}\;\;\Upsilon_{i,j}\\
0&\text{on}\;\;\bdry\Omega_i\smallsetminus\Upsilon_{i,j},
\end{cases}
\end{equation}
and set
\begin{equation}
\label{e:edp-1}
x_{i,n+1}=\lambda_{i,n}x_{i,n}+(1-\lambda_{i,n})(y_{i,n}+a_{i,n}),
\end{equation}
where, for every $i\in\{1,\ldots,m\}$, the following hold.
\begin{enumerate}
\item
\label{c:85}
$x_{i,0}\in \HH_i$.
\item
\label{c:83}
$(a_{i,n})_{n\in\NN}$ and $(b_{i,n})_{n\in\NN}$ are sequences in 
$\HH_i$ such that
\begin{equation}
\label{e:edp+1}
\sum_{n\in\NN}\sqrt{\int_{\Omega_i}|D a_{i,n}|^2}<\pinf
\quad\text{and}\quad 
\sum_{n\in\NN}\sqrt{\int_{\Omega_i}|D b_{i,n}|^2}<\pinf.
\end{equation}
\item
\label{c:84}
$(\lambda_{i,n})_{n\in\NN}$ is a sequence in $\left[0,1\right[$
such that $\sum_{n\in\NN}|\lambda_{i,n}-\lambda_n|<\pinf$.
\end{enumerate}
\end{algorithm}

\subsection{Convergence}

\begin{theorem}
\label{t:s6}
Let $((x_{i,n})_{n\in\NN})_{1\leq i\leq m}$ be sequences generated by
Algorithm~\ref{algo:edp}. Then, for every $i\in\{1,\ldots,m\}$, 
$(x_{i,n})_{n\in\NN}$ converges weakly in $\HH_i$
to a point $x_i\in \HH_i$, and $(x_i)_{1\leq i\leq m}$ is a 
solution to Problem~\ref{prob:edp}. 
\end{theorem}
\begin{proof}
Given $i\in\{1,\ldots,m\}$ and $j\in J(i)$, we first observe that 
${\mathsf T}_{ij}\in\BL(\HH_i,L^2(\Upsilon_{i,j}))$. Indeed, since the 
embedding $\HH_i\hookrightarrow H^1(\Omega_i)$ is continuous 
\cite[p.~1033]{Zeid90} and
${\mathsf T}_{i}\in\BL(H^1(\Omega_i),L^2(\bdry\Omega_i))$, the operator
${\mathsf T}_{ij}\colon\HH_i\to L^2(\Upsilon_{i,j})\colon x\mapsto
({\mathsf T}_ix)|_{\Upsilon_{i,j}}$ is indeed linear and continuous.
Let us now show that Problem~\ref{prob:edp} is a special case of
Problem~\ref{prob:62}. For every $(k,l)\in\KK$ and 
every $i\in\{1,\ldots,m\}$, set
\begin{equation}
\label{e:2009-12-23}
\GG_{kl}=L^2(\Upsilon_{k,l})\quad\text{and}\quad
L_{{kl}i}=
\begin{cases}
{\mathsf T}_{kl},&\text{if}\:\:i=k;\\
-{\mathsf T}_{lk},&\text{if}\:\:i=l;\\
0,&\text{otherwise,}
\end{cases}
\end{equation}
and note that $L_{kli}\in\BL(\HH_i,\GG_{kl})$ since 
\eqref{e:24fevrier2009} entails 
$L^2(\Upsilon_{k,l})=L^2(\Upsilon_{l,k})$. 
Thus, \eqref{prob:16} can be written as
\begin{equation}
\label{e:edpm}
\underset{x_1\in\HH_1,\ldots,x_m\in\HH_m}{\mathrm{minimize}}\;\;
\sum_{i=1}^mf_i(x_i)+\sum_{(k,l)\in\KK}
\varphi_{kl}\bigg(\sum_{i=1}^mL_{{kl}i}x_i\bigg),
\end{equation}
which conforms to \eqref{e:genna07-4}. Next, let us show that 
Algorithm~\ref{algo:edp} is a particular case of Algorithm~\ref{algo:2}. 
To this end, let $i\in\{1,\ldots,m\}$ and $n\in\NN$. Since 
$\bdry\Omega_i=\overline{\Upsilon_{i,i}}\cup
\overline{\bigcup_{j\in J(i)}\Upsilon_{i,j}}$,
we deduce from \cite[Theorem~25.I]{Zeid90} that \eqref{e:edp1} 
admits a unique weak solution $z_{i,n}\in\HH_i$. 
Accordingly \eqref{e:scaledp}, \cite[Definition~25.31]{Zeid90},
\eqref{e:25fevrier2009}, and \eqref{e:L2Ups} yield
\begin{align}
(\forall x\in\HH_i)\quad
\scal{x}{z_{i,n}}
&=\int_{\Omega_i}(Dx)^\top Dz_{i,n}\nonumber\\
&=\int_{\bigcup_{j\in J(i)}\Upsilon_{i,j}}({\mathsf T}_{i}\,x)
\bigg(\sum_{j\in J(i)}\widetilde{v_{ij,n}}
\bigg)\,dS\nonumber\\
&=\sum_{j\in J(i)}\int_{\Upsilon_{i,j}}({\mathsf T}_{ij}\,x)
v_{ij,n}\,dS\nonumber\\
&=\sum_{j\in J(i)}\scal{{\mathsf T}_{ij}\,x}
{v_{ij,n}}_{L^2(\Upsilon_{i,j})}\nonumber\\
&=\Scal{x}{\sum_{j\in J(i)}{\mathsf T}_{ij}^*\,
v_{ij,n}}.
\end{align}
Therefore $z_{i,n}=\sum_{j\in J(i)}{\mathsf T}_{ij}^*v_{ij,n}$ and 
hence \eqref{e:2009-12-23} and \eqref{e:25fevrier2009} yield
\begin{equation}
z_{i,n}=\sum_{k=1}^m\sum_{l\in J(k)}L_{kli}^*\nabla
\varphi_{kl}\bigg(\sum_{j=1}^mL_{klj}x_{j,n}\bigg)
=\sum_{(k,l)\in\KK}L_{kli}^*\nabla
\varphi_{kl}\bigg(\sum_{j=1}^mL_{klj}x_{j,n}\bigg). 
\end{equation}
On the other hand, it follows from \eqref{e:HH_iEDP} and 
\eqref{e:scaledp} that \eqref{e:edp0} is equivalent to
\begin{equation}
\label{e:edp6}
\underset{y\in \HH_i}{\mathrm{minimize}}\;\;
\gamma_nf_i(y)+\frac{1}{2}
\big\|y-\big(x_{i,n}-\gamma_n(z_{i,n}+b_{i,n})\big)\big\|^2,
\end{equation}
the unique solution of which is
\begin{equation}
y_{i,n}=\prox_{\gamma_n f_i}\bigg(x_{i,n}
-\gamma_n\bigg(\sum_{(k,l)\in\KK}L_{kli}^*\nabla\varphi_{kl}
\bigg(\sum_{j=1}^mL_{klj}x_{j,n}\bigg)+b_{i,n}\bigg)\bigg). 
\end{equation}
Moreover, \eqref{e:14nov2007} is implied by \eqref{e:2009-12-23} 
and \eqref{e:edp16oct2008}.
Hence, in view of \eqref{e:edp-1} and \eqref{e:edp+1}, 
Algorithm~\ref{algo:edp} is a particular case of Algorithm~\ref{algo:2}. 
Altogether, Theorem~\ref{t:2} asserts that,
for every $i\in\{1,\ldots,m\}$, the sequence $(x_{i,n})_{n\in\NN}$
converges weakly in $\HH_i$ to a point $x_i\in\HH_i$, where
$(x_i)_{1\leq i\leq m}$ is a solution to Problem~\ref{prob:edp}. 
\end{proof}

\begin{example}
\label{ex:62}
Let $y\in L^2(\Omega)$. With the same notation and hypotheses as in 
Problem~\ref{prob:edp} let, for every $i\in\{1,\ldots,m\}$,
\begin{equation}
f_i\colon\HH_i\to\RR\colon
x\mapsto\frac{1}{2}\int_{\Omega_i}|D x|^2-\int_{\Omega_i}xy
\quad\text{and}\quad(\forall j\in J(i))\quad
\varphi_{ij}=d^2_{C_{ij}},
\end{equation}
where $C_{ij}$ is a nonempty closed convex subset of 
$L^2(\Upsilon_{i,j})$. 
For every $i\in\{1,\ldots,m\}$ the solution to the problem
\begin{equation}
\label{prob:ex2}
\underset{x\in\HH_i}{\mathrm{minimize}}\:f_i(x)
\end{equation}
is the weak solution to the Poisson equation with 
mixed Dirichlet-Neumann conditions \cite[Theorem~25.I]{Zeid90}
\begin{equation}
\begin{cases}
-\Delta x=y&\text{on}\:\:\Omega_i\\
x=0&\text{on}\:\Upsilon_{i,i}\\
\nu_i^\top Dx=0
&\text{on}\;\bigcup_{j\in J(i)}\Upsilon_{i,j}.
\end{cases}
\end{equation}
Problem~\ref{prob:edp} couples these Poisson problems by penalizing 
the violation of the constraints 
${\mathsf T}_{ij}\,x_i-{\mathsf T}_{ji}\,x_j\in C_{ij}$.
\end{example}

\end{document}